\documentclass{elsarticle}
\usepackage{hyperref}
\usepackage[utf8]{inputenc}
\usepackage{wrapfig}
\usepackage{amsmath}
\usepackage{amsthm}
\usepackage{amssymb}
\usepackage[tmargin=1in,bmargin=1in]{geometry}
\usepackage{caption}
\usepackage{subcaption}
\usepackage{graphicx}
\usepackage{latexsym}
\usepackage{pdfsync}
\usepackage[boxed]{algorithm}
\usepackage{algorithm}
\usepackage{multirow}
\usepackage{rotating}
\usepackage{adjustbox}
\usepackage{color}
\usepackage{url}
\usepackage{xcolor,colortbl}
\usepackage{bigstrut}
\usepackage{bigints}
\usepackage[normalem]{ulem}

\usepackage[normalem]{ulem}
\usepackage{booktabs}

\hypersetup{
    colorlinks=true,              
    linkcolor=blue,               
    citecolor=blue,               
    filecolor=black,              
    urlcolor=red,               
    bookmarks=false,
    pdffitwindow=true,
    pdfpagelayout=SinglePage
}

\definecolor{mypink}{cmyk}{0, 0.7808, 0.4429, 0.1412}

\definecolor{mybrown}{cmyk}{0, 0.60, 0.95, 0.63}
\definecolor{darkgreen}{cmyk}{0.98, 0, 0.36, 0.22}

\newcommand{\be}{\begin{equation}}
\newcommand{\ee}{\end{equation}}
\newcommand{\ben}{\begin{equation*}}
\newcommand{\een}{\end{equation*}}
\newcommand{\bea}{\begin{eqnarray}}
\newcommand{\eea}{\end{eqnarray}}
\newcommand{\bean}{\begin{eqnarray*}}
\newcommand{\eean}{\end{eqnarray*}}

\usepackage{algorithm}
\usepackage[noend]{algpseudocode}
\makeatletter
\def\BState{\State\hskip-\ALG@thistlm}
\makeatother

\begin{document}

\title{JAX-DIPS: Neural bootstrapping of finite discretization methods and application to elliptic problems with discontinuities}

\cortext[cor]{Corresponding author: p.a.mistani@gmail.com}

\author[1,2]{Pouria A. Mistani\thanks{corresponding author} $^{\dagger,}$}
\author[1]{Samira Pakravan$^{\dagger,}$}
\author[2]{Rajesh Ilango}
\author[1]{Frederic Gibou}

\address[1]{University of California, Santa Barbara, CA 93106, USA}
\address[2]{NVIDIA Corporation, Santa Clara, CA 95051, USA}

\begin{abstract}
	We present a scalable strategy for development of mesh-free hybrid neuro-symbolic partial differential equation solvers based on existing mesh-based numerical discretization methods. Particularly, this strategy can be used to efficiently train neural network surrogate models of partial differential equations by (i) leveraging the accuracy and convergence properties of advanced numerical methods, solvers, and preconditioners, as well as (ii) better scalability to higher order PDEs by strictly limiting optimization to first order automatic differentiation. The presented neural bootstrapping method (hereby dubbed NBM) is based on evaluation of the finite discretization residuals of the PDE system obtained on implicit Cartesian cells centered on a set of random collocation points with respect to trainable parameters of the neural network. Importantly, the conservation laws and symmetries present in the bootstrapped finite discretization equations inform the neural network about solution regularities within local neighborhoods of training points. We apply NBM to the important class of elliptic problems with jump conditions across irregular interfaces in three spatial dimensions. We show the method is convergent such that model accuracy improves by increasing number of collocation points in the domain and predonditioning the residuals. We show NBM is competitive in terms of memory and training speed with other PINN-type frameworks. The algorithms presented here are implemented using \texttt{JAX} in a software package named \texttt{JAX-DIPS} (\href{https://github.com/JAX-DIPS/JAX-DIPS}{https://github.com/JAX-DIPS/JAX-DIPS}), standing for differentiable interfacial PDE solver. We open sourced \texttt{JAX-DIPS} to facilitate research into use of differentiable algorithms for developing hybrid PDE solvers.

\end{abstract}


\maketitle
\def\thefootnote{$\dagger$}\footnotetext{These authors contributed equally to this work.}

\section{Introduction}
\label{sec::introduction}

\subsection{Problem statement}
Consider a closed irregular interface ($\rm \Gamma$) that partitions the computational domain ($\rm \Omega$) into interior ($\rm \Omega^-$) and exterior ($\rm \Omega^+$) subdomains; \textit{i.e.}, $\rm \Omega=\Omega^- \cup \Gamma \cup \Omega^+$. We are interested in the solutions $\rm u^\pm\in \Omega^\pm$ to the following class of linear elliptic problems in  $\rm \mathbf{x}\in\Omega^\pm$:
\begin{align*}
	 & k^{\pm}u^{\pm} - \nabla \cdot (\mu^{\pm}\nabla u^\pm)=f^{\pm}, & \mathbf{x}\in\Omega^\pm \\
	 & [u]=\alpha,                                                    & \mathbf{x} \in \Gamma   \\
	 & [\mu \partial_{\mathbf{n}}u]=\beta,                            & \mathbf{x} \in \Gamma
\end{align*}
Here $f^\pm=f(\mathbf{x} \in \Omega^\pm)$ is the spatially varying source term, $\rm \mu^\pm=\mu(\mathbf{x} \in \Omega^\pm)$  are the diffusion coefficients, and $k^\pm$ are the reaction coefficients in the two domains. We consider Dirichlet boundary conditions in a cubic domain $\rm \Omega=[-L/2,L/2]^3$.

This class of problems arise ubiquitously in describing diffusion dominated processes in physical systems and life sciences where sharp and irregular interfaces regulate transport across regions with different properties. Examples include Poisson-Boltzmann equation for describing electrostatic properties of membranes, colloids and solvated biomolecules with jump in dielectric permittivities \cite{sharp1990calculating,MirzadehPB}, electroporation of cell aggregates with nonlinear membrane jump conditions \cite{mistani2019parallel}, epitaxial growth in fabrication of opto-electronic devices where atomic islands grow by surface diffusion of adatoms across freely moving interfaces \cite{MISTANI2018150}, solidification of multicomponent alloys used for manufacturing processes with free interfaces separating different phases of matter \citep{theillard2015sharp,bochkov2021sharp}, directed self-assembly of diblock copolymers for next generation lithography \cite{galatsis2010patterning,ouaknin2018level,bochkov2021non}, multiphase flows with and without phase change, and Poisson-Nernst-Planck equations for electrokinetics. Much of these processes are multiscale and the changes across interfaces must be mathematically modeled and numerically solved as sharp surfaces. Smoothing strategies introduce unphysical characteristics in the solution and lead to systemic errors.

\subsection{Literature on relevant finite discretization methods}
Several numerical methods have been proposed for accurate solution of this class of problems based on explicit or implicit representation of the interface. Finite element methods rely on explicit meshing of the surface that poses severe challenges \cite{babuvska1970finite,bramble1996finite}. Implicit methods include the Immersed Interface Method (IIM) \cite{leveque1994immersed} and its variants \cite{adams2002immersed,li2003new,ewing1999immersed,gong2008immersed} that rely on Taylor expansions of the solution on both sides of the interface and modifying the local stencils to impose the jump conditions. The main challenge is evaluating high order jump conditions and surface derivatives along interface. Another method is the Ghost Fluid Method (GFM) \cite{fedkiw1999non} that was originally introduced to approximate two-phase compressible flows and later applied to the Poisson problem with jump conditions \cite{liu2000boundary}. The basic idea is to define fictitious fluid regions across the discontinuities by adding jump conditions to the true fluid. While GFM captures the normal jump in solution accurately, the tangential jump is smeared. This was solved by the Voronoi Interface Method (VIM) \cite{guittet2015solving} by applying the GFM treatment on a local Voronoi mesh by adapting a local Cartesian mesh which introduces numerical challenges. Several other approaches include the cut-cell method \cite{crockett2011cartesian}, discontinuous Galerkin and eXtended Finite Element Method (XFEM) \cite{lew2008discontinuous,moes1999finite,belytschko2001arbitrary} among others.

In this work we bootstrap the level-set based finite volume method on Cartesian grids proposed by Bochkov \& Gibou (2020) \cite{BOCHKOV2020109269}. This method is based on the idea of Taylor expansions in the normal direction and employing one-sided least-square interpolations for imposing jump conditions. In particular, this method offers second order accurate numerical solutions with first order accurate gradients in the $\rm L^\infty$-norm.

\subsection{Literature on solving PDEs with neural networks}
Since early 1990s, artificial neural networks have been used for solving partial differential equations by (i) mapping the algebraic operations of the discretized PDE systems onto specialized neural network architectures and minimizing the network energy, or (ii) treating the whole neural network as the basic approximation unit whose parameters are adjusted to minimize a specialized error function that includes the differential equation itself with its boundary/initial conditions.

In the first category, neurons output the discretized solution values over a set number of grid points and minimizing the network energy drives the neuronal values towards the solution of the linear system at the mesh points. In this case, the neural network energy is the residual of the finite discretization method summed over all neurons of the network \cite{lee1990neural}. Although the convergence properties of the finite discretization methods gaurantee and control quality of the obtained solutions, the computational costs grow by increasing resolution and dimensionality. Interestingly, due to regular and sparse structure of the finite discretizations, such locally connected neural network PDE solvers have been implemented on VLSI analog CMOS circuits \cite{gobovic1993design,chua1988cellular,chua1988cellularA}.

The second strategy proposed by Lagaris \textit{et al.} \cite{lagaris1998artificial} relies on the function approximation capabilities of the neural networks. Encoding the solution everywhere in the domain within a neural network offers a mesh-free, compact, and memory efficient surrogate model for the solution function that can be utilized in subsequent inference tasks. This method has recently re-emerged as the physics-informed neural networks (PINNs) \cite{RAISSI2019686} and is widely used. 

Despite their implementation simplicity and offering fast inference on accelerated hardware, these methods suffer from several shortcomings:

\begin{enumerate}
\item lack controllable accuracy and convergence properties of finite discretization methods \cite{krishnapriyan2021characterizing}, 

\item computing the loss and optimizing it involves evaluation of second order (and higher order) gradients using automatic differentiation (AD) through deep neural networks which leads to evaluating exponentially large computational graphs that is extremely memory-intensive, slow, and impractical to scale,

\item the basic assumption that automatic differentiation capabilities of current machine learning frameworks can evaluate ``exact'' derivatives across complex surrogate models is fundamentally flawed \cite{johnson2023software},

\item automatic differentiation of an \textit{un-optimized} neural network during training to compute the spatial gradients does not offer exact gradients for the PDE. Although these derivatives are ``exact'' (see \cite{johnson2023software} for a discussion) within the parameters of the neural network , it is important to note that these derivatives do not represent the true spatial derivatives of the solution but \textit{exact derivatives of an approximate function}.

\end{enumerate}

These shortcomings motivate pursuit of hybrid solvers to combine the performance advantages of neural network inference on modern accelerated hardware with the accuracy of finite discretization methods while reducing the computational costs and errors associated with excessive use of AD. 

Hybridization efforts are algorithmic or architectural. One important algorithmic method is the deep Galerkin method (DGM) \cite{SIRIGNANO20181339} that is a neural network extension of the mesh-free Galerkin method where the solution is represented as a deep neural network rather than a linear combination of basis functions. The mesh-free nature of DGM, that stems from the underlying mesh-free Galerkin method, enables solving problems in higher dimensions by training the neural network model to satisfy the PDE operator and its initial and boundary conditions on a randomly sampled set of points rather than on an exponentially large grid. Although the number of points is huge in higher dimensions, the algorithm can process training on smaller batches of data points sequentially. Besides, second order derivatives in PDEs are calculated by a Monte Carlo method that retain scaling to higher dimensions. Another important algorithmic method is the deep Ritz method for solving variational problems \cite{yu2018deep} that implements a deep neural network approximation of the trial function that is constrained by numerical quadrature rule for the variational functional, followed by stochastic gradient descent.

Architectural hybridization methods are based on differentiable numerical linear algebra. One emerging class involves implementing differentiable finite discretization solvers and embedding them in the neural network architectures that enable application of end-to-end differentiable gradient based optimization methods. Recently, differentiable solvers have been developed in JAX \cite{jax2018github} for fluid dynamic problems, such as \texttt{Phi-Flow} \cite{holl2020phiflow}, \texttt{JAX-CFD} \cite{Kochkov2021-ML-CFD}, and \texttt{JAX-FLUIDS} \cite{bezgin2022jax}. These methods are suitable for inverse problems where an unknown field is modeled by the neural network, while the model influence is propagated by the differentiable solver into a measurable residual \cite{pakravan2021solving,dal2020data,lu2020extracting}. We also note the classic strategy for solving inverse problems is the adjoint method to obtain the gradient of the loss without differentiation across the solver \cite{berg2017neural}; however, deriving analytic expression for the adjoint equations is tedious, should be repeated after modification of the problem or its loss function, and can become impractical for multiphysics problems. Other important utilities of differentiable solvers are to model and correct for the solution errors of finite discretization methods \cite{um2020solver}, learning and controling PDE systems \cite{de2018end,holl2020learning}.

Neural networks are not only universal approximators of continuous functions, but also of nonlinear operators \cite{chen1995universal}. Although this fact has been leveraged using data-driven strategies for learning differential operators by many authors \cite{lu2019deeponet,bhattacharya2020model,li2020neural,li2020fourier}, current authors have demonstrated utility of differentiable solvers to effectively train nonlinear operators without any data in a completely physics-driven fashion, see section on learning the inverse transforms in \cite{pakravan2021solving}. In subsequent work we will demonstrate how NBM can be used to train neural operators in a purely physics-driven fashion.

In this work we propose a novel algorithm for solving PDEs based on deep neural networks by lifting any existing mesh-based finite discretization method off of its underlying grid and extend it into a mesh-free method that can be applied to high dimensional problems on unstructured random points in an embarrasingly parallel fashion. In section \ref{sec:nbm} we present the neural bootstrapping method, next we apply it to an advanced finite volume discretization scheme for elliptic problems with jump conditions across irregular geometries in section \ref{sec:dips}. We show numerical results of the proposed framework on interfacial PDE problems in section \ref{sec:examples} and conclude with section \ref{sec:conc}.

\section{Neural Bootstrapping Method (NBM)} \label{sec:nbm}

\subsection{Algorithm}
\begin{figure}
     \centering
     \begin{subfigure}[b]{\textwidth}
         \centering
         \includegraphics[width=\textwidth]{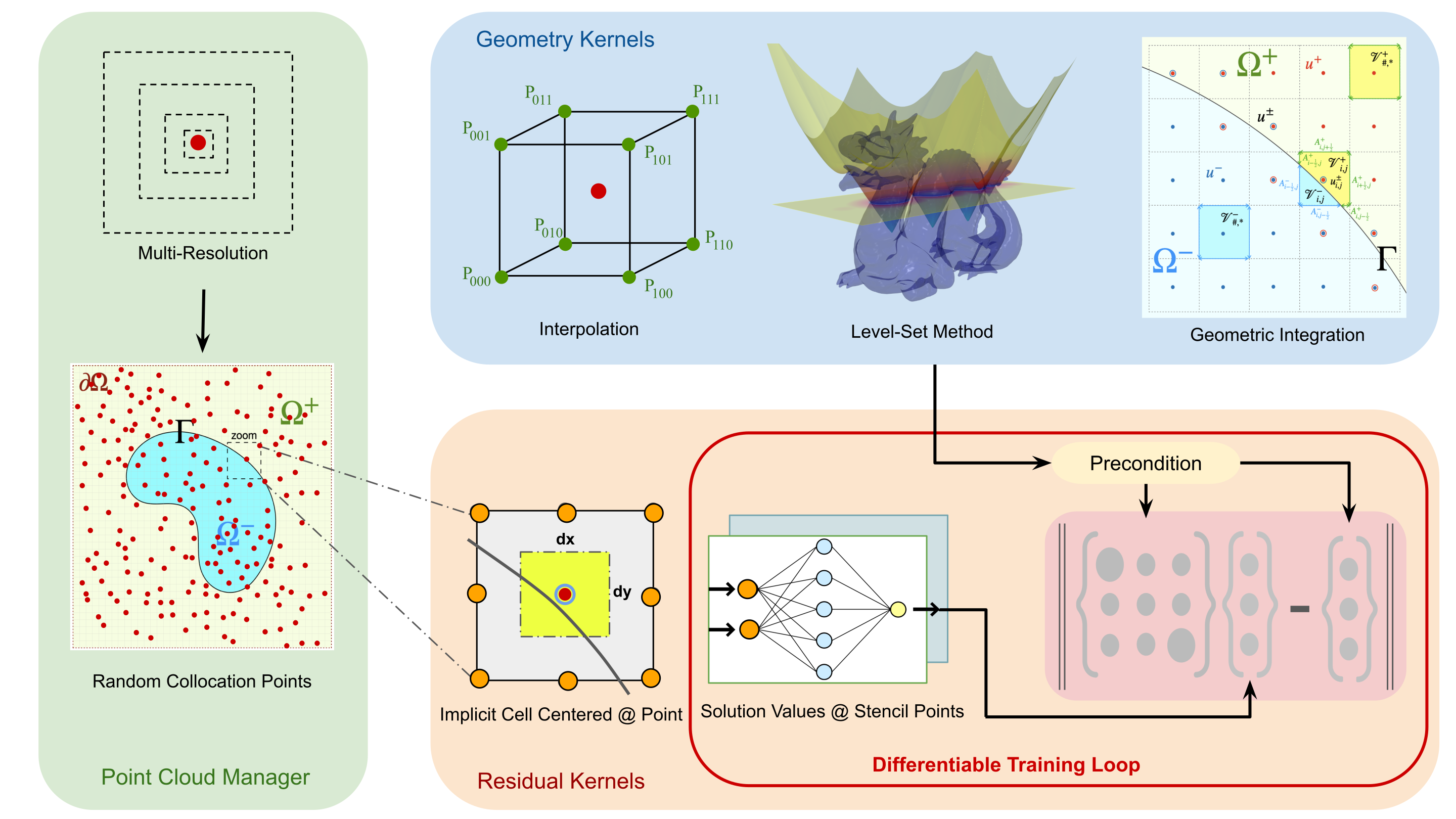}
         \caption{NBM kernels compute residual contribution by each collocation point per thread. Kernel operations involve considering implicit cells at different resolutions according to the bootstrapped finite discretization method. The point-wise evaluations at each implicit cell is locally preconditioned based on the geometry of the interfaces crossing the implicit cells.}
         \label{fig:nbm_kernel}
     \end{subfigure}
     \hfill
     \begin{subfigure}[b]{\textwidth}
         \centering
         \includegraphics[width=\textwidth]{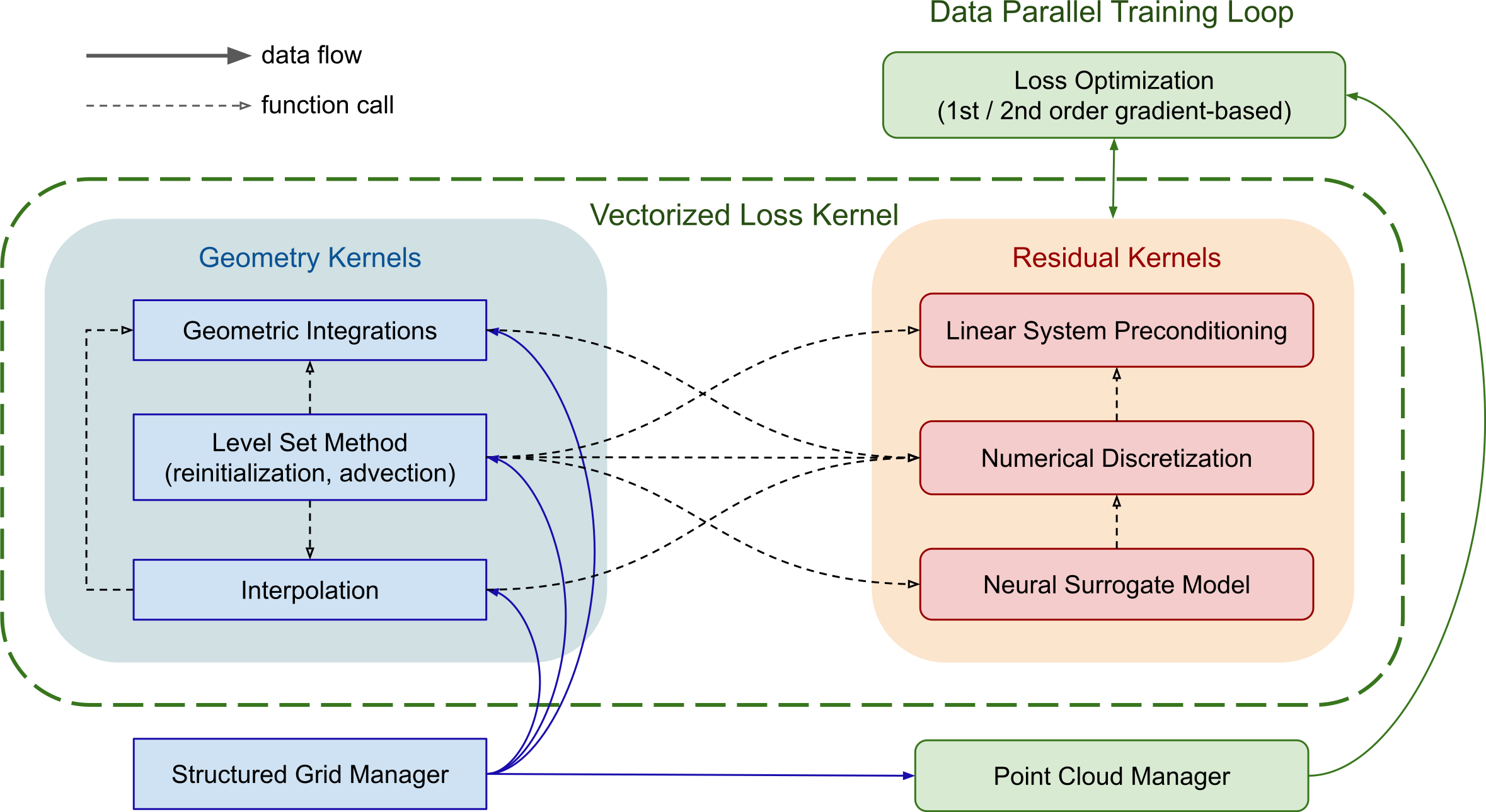}
         \caption{\texttt{JAX-DIPS} software architecture layout. Geometric information is managed by a mesh oracle that is a structured mesh at much lower resolution that stores the level-set function. The training loop involves automatic differentiation across the assembly of the linear system performed by the NBM kernels. Data distribution is achieved by composing the point-wise loss kernel with \texttt{jax.pmap} and \texttt{jax.vmap}.}
         \label{fig:jax_dips}
     \end{subfigure}
        \caption{Neural Bootstrapping Method (NBM) and the \texttt{JAX-DIPS} software architecture.}
        \label{fig:nbm}
\end{figure}

Figure \ref{fig:nbm} illustrates schematic of the proposed algorithm. Neural networks are used as surrogates for the solution function that are iteratively adjusted to minimize discretization residuals at a set of randomly sampled points and at arbitrary resolutions. The key idea is that neural networks can be evaluated over vertices of any discretization stencil centered at any point in the domain to emulate the effect of an structured mesh without ever materializing the mesh. Therefore, we use neural networks to bootstrap mesh based finite discretization (FD) methods to compile mesh free numerical methods. Operations in differentiable NBM kernels are:
\begin{enumerate}

\item A compute cell is implicitly constructed at any input coordinate and at a specified resolution. At the presence of discontinuities a coarse mesh encapsulates an interpolant for the level-set function whose intersection with the implicit cell is calculated to obtain necessary geometric information for the FD kernel and preconditioner.

\item FD kernel is applied on the compute cell where the solution values are evaluated by the neural network. Each kernel contributes a local linear system $\rm L^2$-norm residual $\rm r_p = ||A u_p - b||$ at a point $\rm p$.

\item Residuals are preconditioned using common preconditioners to balance relative magnitude of contributions from different points and set them on equal level of importance before summing to produce a global loss value.

\item Gradient based optimization methods used in machine learning are applied to adjust neural network parameters. The automatic differentiation loop passes across the NBM kernels, see figure \ref{fig:nbm_kernel}.

\end{enumerate}

NBM training of neural network surrogate models for PDEs offers several benefits:

\begin{itemize}
\item FD methods offer gauranteed accuracy and controllable convergence properties for the training of neural network surrogate models. These are critical features for solving real-world complex physical systems using neural networks.

\item NBM offers a straightforward path for applying mesh-based FD methods on unstructured random points. This is an important ability for augmenting observational data in the training pipelines.

\item The algorithm is highly parallelizable and is ideally suited for GPU-accelerated computing paradigm. 

\item Multi-GPU parallel solution of PDE systems is reduced to the much simpler problem of data-parallel training using existing machine learning frameworks. Data parallelism involves distributing collocation points across multiple processors to compute gradient updates and then aggregating these locally computed updates \cite{shallue2018measuring}.

\item Only first order automatic differentiation is required for training PDE systems. This dramatically reduces memory requirements and computational costs associated with higher order AD computations across neural network models in other methods.

\item Use of first order optimizers, enabled by differentiable finite discretizations, could improve scaling of traditional PDE solvers that use second order optimization techniques.

\end{itemize}

\subsection{Neural network approximators for the solution}

In 1987, Hecht and Nielson \cite{hecht1987kolmogorov} applied an improved version of Kolmogorov's 1957 superposition theorem \cite{kolmogorov1957representation}, due to Sprecher \cite{sprecher1965structure}, to the field of neurocomputing and demonstrated that a $3$-layer feedforward neural network (one input layer with $n$ inputs, one hidden layer with $2n+1$ neurons, one output layer) are universal approximators for all \textit{continuous} functions from the $n$-dimensional cube to a finite $m$-dimensional real vector space; \textit{i.e.}, $f: [0,1]^n \rightarrow \mathbb{R}^m $. Recently, Ismailov (2022) \cite{ismailov2022} demonstrated existence of neural networks implementing \textit{discontinuous} functions, however efficient learning algorithms for such networks are not still available. 

The solutions of interfacial PDE problems are discontinuous, with jumps appearing not only in the solution but also in the solution gradient. In light of above considerations, we define two separate neural networks to represent solution in $\Omega^-$ and $\Omega^+$ regions:
\begin{align*}
& u^+ = \mathcal{N}^+(\mathbf{x}): \mathbb{R}^3\cap \Omega^+ \rightarrow \mathbb{R}  & u^- = \mathcal{N}^-(\mathbf{x}): \mathbb{R}^3\cap \Omega^- \rightarrow \mathbb{R}
\end{align*}

We use SIREN neural networks, where we implement fully connected feedforward architecture with \texttt{sin} activation function and the output layer is a single linear neuron. Note that piecewise differentiable nonlinearities such as the \texttt{ReLU} function are inappropriate choices for representing solutions to differential equations. Weights and biases are initialized from a truncated normal distribution with zero mean and unit variance. 

\begin{figure}[ht]
	\centering
	\includegraphics[width=\linewidth]{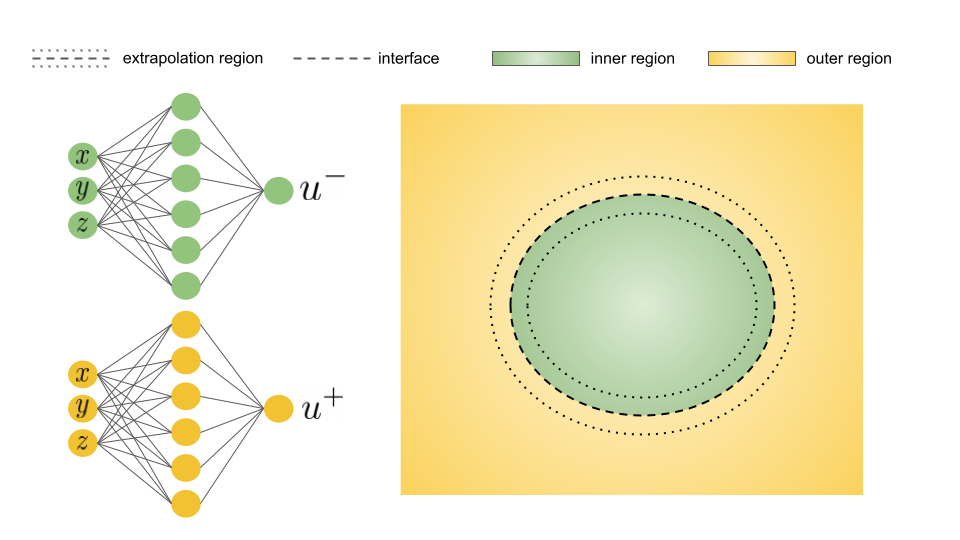}
	\caption{Two neural networks are defined for the two regions of the computational domain.}
	\label{fig:shapes}
\end{figure}

Solution networks are evaluated on sampled points in the domain while the parameters of these networks are optimized using the loss function. We define the loss function by the mean-squared-error (MSE) of the residual of the discretized partial differential equation with jump conditions derived in section \ref{sec::FD} that is evaluated on the grid points:
\begin{align*}
	&\mathcal{L}(u) = \vert\vert A\hat{u}_{\theta}(\mathbf{x}_{ijk}\in \Omega) - b\vert\vert_2^2  
\end{align*}
However, other choices such as \texttt{Huber} \cite{huber1992robust} or \texttt{log-cosh} \cite{saleh2022statistical} loss functions may improve results by automatically suppressing $\rm L_1$ norm for larger residual values while minimizing $\rm L_2$ norm for smaller values of the residual. \texttt{JAX-DIPS} allows for computation of the gradient of the loss function using automatic differentiation, \textit{i.e.} $\nabla_\theta \mathcal{L}(u)$ where $\theta$'s are network parameters. Therefore, our strategy is to leverage this capability and use first order optimizers developed in the deep learning community (such as Adam \cite{kingma2014adam}, \textit{etc}) to minimize the aforementioned loss function. We emphasize the main benefit of using first order gradient based optimization algorithms is better memory efficiency that is suitable for large scale optimization problems with very large number of parameters in the neural network model.

In the remainder of this manuscript we present details of applying NBM to solving elliptic problems with discontinuities across irregular interfaces.
\section{JAX-DIPS: Differentiable Interfacial PDE Solver} \label{sec:dips}
We developed an end-to-end differentiable library for solving the elliptic problems with discontinuities in solution and solution gradient across irregular geometries. In \texttt{JAX-DIPS}, we bootstrap a sophisticated and modern finite volume discretization method \cite{BOCHKOV2020109269}. The geometries are represented implicitly using the level-set function on a coarse mesh. We have implemented a uniform grid that supports operations such as interpolations, interface advection, integrations over interfaces as well as in domains. We describe the numerical algorithms for the level-set module and the elliptic solver in this section. Figure \ref{fig:jax_dips} illustrates a high-level overview of the \texttt{JAX-DIPS} software architecture.

Below we present and compare two possible approaches for treating the jump conditions in the interfacial PDE solver: (i) regression-based extrapolation, and (ii) neural extrapolation. The main difference between the two approaches is that approach (i) only requires first order AD for optimizing the loss, while (ii) effectively requires second order AD computations due to first evaluation of the loss and a second AD during optimization of the loss. Approach (i) offers better computational properties thanks to a complete bootstrapping of the underlying finite volume discretization method.

\subsection{Approach I. Finite discretization method fused with regression-based extrapolation}\label{sec::FD}

For spatial discretizations at the presence of jump conditions we employ the numerical algorithm proposed by Bochkov and Gibou (2020) \cite{bochkov2020solving}. This method produces second-order accurate solutions and first-order accuracte gradients in the $L^\infty$-norm, while having a compact stencil that makes it a good candidate for parallelization. Moreover, treatment of the interface jump conditions do not introduce any augmented variables, this preserves the homogeneous structure of the linear system. Most importantly, jump conditions only appear on the right-hand-side of the discretization and do not pollute the matrix term, this is beneficial for accelerating the solver. Here we use a background uniform 2D grid only for presentation of the finite volume discretization equations; we will not use this grid in the actual implementation but instead assume a local 3D cell around random points spanning in the domain during optimization. 

At points where the finite volumes are crossed by $\rm \Gamma$ we have
\begin{align*}
	 & \sum_{s=-,+}\int_{\Omega^s \cap \mathcal{V}_{i,j}}  k^{s} u^{s} d\Omega -\sum_{s=-,+} \int_{\Omega^s\cap \partial\mathcal{V}_{i,j}} \mu^{s}\partial_{\mathbf{n}^s} u^s  d\Gamma = \sum_{s=-,+}\int_{\Omega^s \cap \mathcal{V}_{i,j}}  f^{s} d\Omega + \int_{\Gamma \cap \mathcal{V}_{i,j}}[\mu\partial_{\mathbf{n}}u]d\Gamma \\
	\intertext{following standard treatment of volumetric integrals and using central differencing for derivatives we obtain in 2D (with trivial 3D extension)}
	 & \sum_{s=-,+} k_{i,j}^s u_{i,j}^{s} |\mathcal{V}_{i,j}^s| - \sum_{s=-,+}\bigg( \mu_{i-\frac{1}{2},j}^s A_{i-\frac{1}{2},j}^s\frac{u_{i-1,j}^s - u_{i,j}^s}{\Delta x}     +   \mu_{i+\frac{1}{2},j}^s A_{i+\frac{1}{2},j}^s\frac{u_{i+1,j}^s - u_{i,j}^s}{\Delta x} +                                                          \\
	 & \mu_{i, j-\frac{1}{2}}^s A_{i, j-\frac{1}{2}}^s\frac{u_{i,j-1}^s - u_{i,j}^s}{\Delta y} + \mu_{i, j+\frac{1}{2}}^s A_{i, j+\frac{1}{2}}^s\frac{u_{i,j+1}^s - u_{i,j}^s}{\Delta y} \bigg)                                                                                                                                     \\
	 & =  \sum_{s=-,+} f_{i,j}^{s} |\mathcal{V}_{i,j}^s| + \int_{\Gamma\cap \mathcal{V}_{i,j}} \beta d\Gamma + \mathcal{O}(\max(\Delta x, \Delta y)^\mathcal{D})
\end{align*}
where $\mathcal{D}$ is the problem dimensionality.
\begin{figure}
	\centering
	\includegraphics[width=0.45\linewidth]{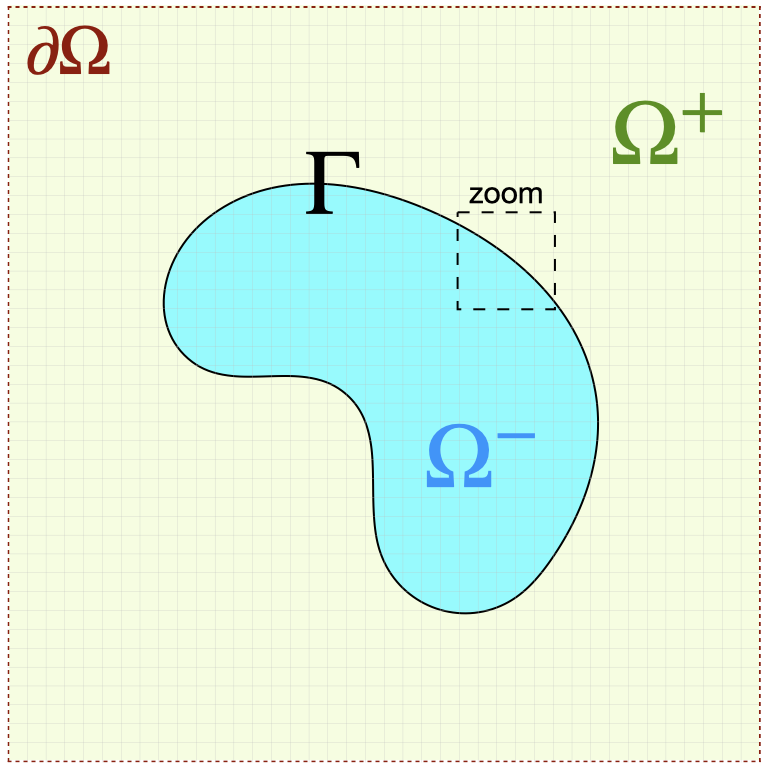}
	\includegraphics[width=0.45\linewidth]{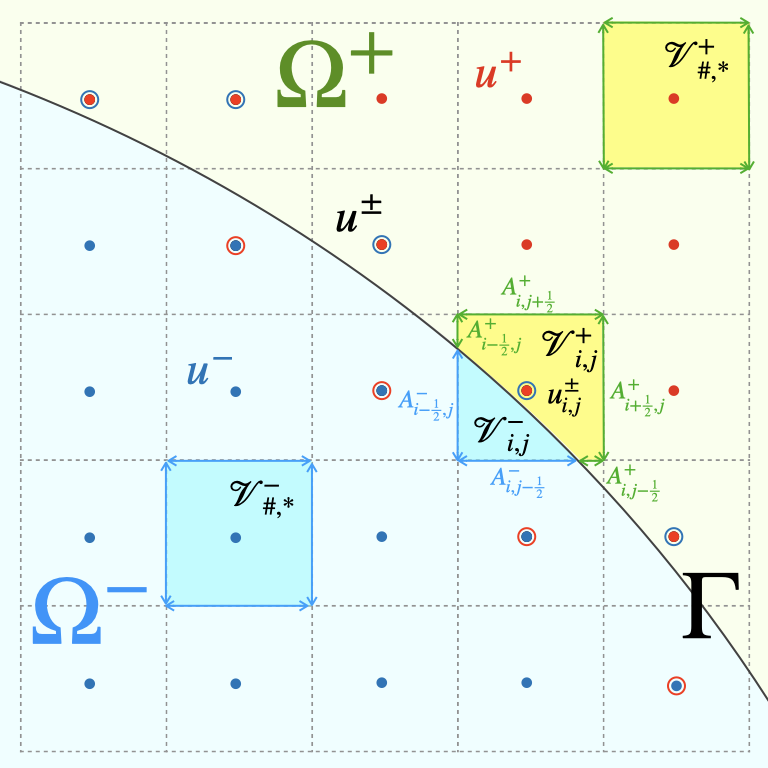}
	\caption{Notation used in this paper. Close to the interface where finite volumes are crossed by the interface, there are extra degrees of freedom (open circles) that are extrapolations of solutions from each domain to the opposite domain. Jump conditions are implicitly encoded in these extrapolated values.}
	\label{fig:grid}
\end{figure}
Note that far from interface either $s=-$ (for $\mathbf{x}\in \Omega^-$) or $s=+$ (for $\mathbf{x}\in \Omega^+$) is retained. This is automatically considered through zero values for sub-volumes $|\mathcal{V}_{i,j}^+|$ and $|\mathcal{V}_{i,j}^-|$ as well as their face areas. Note that $\mu_{i-1/2,j}^-$ (or $\mu_{i-1/2,j}^+$) corresponds to the value of diffusion coefficient at the middle of segment $A^-_{i-1/2,j}$ (or $A^+_{i-1/2,j}$) respectively, same is true for other edges as well. However, there are extra degrees of freedom on grid points whose finite volumes are crossed by the interface; \textit{i.e.}, see double circles in figure \ref{fig:grid}. \cite{bochkov2020solving} derived analytical expressions for the extra degrees of freedom ($u^+$ in $\Omega^-$ and $u^-$ in $\Omega^+$) in terms of the original degrees of freedom ($u^-$ in $\Omega^-$ and $u^+$ in $\Omega^+$) as well as the jump conditions, this preserves the original $\rm N_x \times N_y$ system size. The basic idea is to extrapolate the jump at grid point from jump condition at the projected point onto the interface using a Taylor expansion: $u^+_{i,j} - u^-_{i,j}=[u]_{\mathbf{r}_{i,j}^{pr}} + \delta_{i,j}(\partial_\mathbf{n}u^+(\mathbf{r}^{pr}_{i,j}) - \partial_\mathbf{n}u^-(\mathbf{r}^{pr}_{i,j})) $. The unknown value ($u^-_{i,j}$ or $u^+_{i,j}$) is obtained based on approximation of the normal derivatives (\textit{i.e.} $\partial_\mathbf{n}u^\pm(\mathbf{r}^{pr}_{i,j})$) which are computed using a least squares calculation on neighboring grid points that are in the fast-diffusion region (referred to as ``Bias Fast'') or in the slow diffusion region (referred to as ``Bias Slow''). This makes two sets of rules for unknown values $u^\pm_{i,j}$.

In two dimensions and on uniform grids, the gradient operator at the grid cell $(i,j)$ that is crossed by an interface is estimated by a least squares solution given by
\begin{align*}
	 & (\nabla u^\pm)_{i,j} = \mathbf{D}^\pm_{i,j} \begin{bmatrix}
		u_{i-1,j-1} - u^\pm_{i,j} \\
		u_{i,j-1} - u^\pm_{i,j}   \\
		\vdots                    \\
		u_{i+1,j+1} - u^\pm_{i,j}
	\end{bmatrix} & \mathbf{D}^\pm_{i,j} = \big(X^T_{i,j} W^\pm_{i,j} X_{i,j} \big)^{-1} \big( W^\pm_{i,j} X_{i,j} \big)^T
\end{align*}
and
\begin{align*}
	 & W^\pm_{i,j} = \begin{bmatrix}
		 & \omega^\pm_{i,j} (-1,-1) &                         &        &                        \\
		 &                          & \omega^\pm_{i,j} (0,-1) &        &                        \\
		 &                          &                         & \ddots &                        \\
		 &                          &                         &        & \omega^\pm_{i,j} (1,1) \\
	\end{bmatrix} & X_{i,j} = \begin{bmatrix}
		 & -h_x & -h_y \\
		 & 0    & -h_y \\
		 & h_x  & -h_y \\
		 & -h_x & 0    \\
		 & 0    & 0    \\
		 & h_x  & 0    \\
		 & -h_x & h_y  \\
		 & 0    & h_y  \\
		 & h_x  & h_y
	\end{bmatrix}
\end{align*}
and
\begin{equation}
	\omega_{i,j}^\pm (p,q) = \begin{cases}
		1 & (p,q)\in N_{i,j}^\pm \\
		0 & \text{else}
	\end{cases}
\end{equation}
In this case, $D^\pm_{i,j}$ is a $2\times 9$ matrix and we denote each of its $2\times 1$ columns with $d^\pm_{i,j,p,q}$
\begin{align*}
	\mathbf{D}^\pm_{i,j}  = \begin{bmatrix}
		 & d^\pm_{i,j,-1,-1} & d^\pm_{i,j,0,-1} & d^\pm_{i,j,1,-1} & d^\pm_{i,j,-1,0} & d^\pm_{i,j,0,0} & d^\pm_{i,j,1,0} & d^\pm_{i,j,-1,1} & d^\pm_{i,j,0,1} & d^\pm_{i,j,1,1}
	\end{bmatrix}
\end{align*}
The least square coefficients are then obtained by dot product of normal vector with these columns
\begin{align*}
	c^\pm_{i,j,p,q} = \mathbf{n}_{i,j}^T  d^\pm_{i,j,p,q}
\end{align*}
and normal derivative can be computed (noting that $c_{i,j}^\pm =-\sum_{(p,q)\in N_{i,j}^\pm}c_{i,j,p,q}^\pm$)
\begin{align*}
\partial_n u^\pm (\mathbf{r}^{proj}_{i,j})= c_{i,j}^\pm u_{i,j}^\pm + \sum_{(p,q)\in N_{i,j}^\pm} c_{i,j,p,q}^\pm u_{i+p, j+q}^\pm + \mathcal{O}(h) 
\end{align*}

At this point we can define a few intermediate variables at each grid point to simplify the presentation of the method,
\begin{align*}
	 & \zeta_{i,j,p,q}^\pm := \delta_{i,j} \frac{[\mu]}{\mu^\mp}c_{i,j,p,q}^\pm  & \zeta_{i,j}^\pm := -\sum_{(p,q)\in N_{i,j}^\pm} \zeta_{i,j,p,q}^\pm   \\
	 & \gamma_{i,j,p,q}^\pm := \frac{\zeta_{i,j,p,q}^\pm}{1 \pm \zeta^\pm_{i,j}} & \gamma^\pm_{i,j} := -\sum_{(p,q)\in N_{i,j}^\pm} \gamma_{i,j,p,q}^\pm
\end{align*}
where the set of neighboring grid points are
\begin{align*}
	N_{i,j}^\pm = \{(p,q) :\ \ \  p=-1,0,1, \ \ \  q=-1,0,1, \ \ \ (p,q)\neq (0,0), \ \ \ \mathbf{x}_{i+p,j+q}\in \Omega^\pm \}
\end{align*}
and $\delta_{i,j}$ is the signed distance from $\mathbf{x}_{i,j}$ that is computed from the level-set function $\phi(\mathbf{x})$
\begin{align*}
	\delta_{i,j}=\frac{\phi(\mathbf{x}_{i,j})}{|\nabla \phi(\mathbf{x}_{i,j})|}
\end{align*}

\begin{itemize}
	\item Rules based on approximating $\partial_\mathbf{n}u^+(\mathbf{r}^{pr}_{i,j})$:
\end{itemize}

\begin{equation}
	u_{i,j}^-=\begin{cases}
		u_{i,j}                                                                                                                                                 & \text{ $\mathbf{x}_{i,j}\in \Omega^-$} \\
		u_{i,j} (1 - \gamma_{i,j}^- ) - \sum_{(p,q)\in N_{i,j}^-} \gamma_{i,j,p,q}^- u_{i+p,j+q} - (\alpha + \frac{\delta_{i,j}\beta}{\mu^+})(1-\gamma^-_{i,j}) & \text{ $\mathbf{x}_{i,j}\in \Omega^+$}
	\end{cases}
\end{equation}
\begin{equation}
	u_{i,j}^+=\begin{cases}
		u_{i,j}(1 - \zeta^-_{i,j} ) - \sum_{(p,q)\in N_{i,j}^-} \zeta^-_{i,j,p,q}u_{i+p,j+q}  + \alpha + \delta_{i,j}\frac{\beta}{\mu^+} & \text{ $\mathbf{x}_{i,j}\in \Omega^-$} \\
		u_{i,j}                                                                                                                          & \text{ $\mathbf{x}_{i,j}\in \Omega^+$}
	\end{cases}
\end{equation}
It is useful to cast this in the form of matrix kernel operations through defining intermediate tensors:
\begin{align*}
	 & \boldsymbol{\Gamma}_{i,j} := \begin{bmatrix}
		 & \gamma_{i-1,j+1}^- & \gamma^-_{i,j+1} & \gamma^-_{i+1,j+1} \\
		 & \gamma_{i-1,j}^-   & \gamma^-_{i,j}   & \gamma^-_{i+1,j}   \\
		 & \gamma_{i-1,j-1}^- & \gamma^-_{i,j-1} & \gamma^-_{i+1,j-1}
	\end{bmatrix}, & \boldsymbol{\zeta}_{i,j} := \begin{bmatrix}
		 & \zeta^-_{i-1,j+1} & \zeta^-_{i,j+1} & \zeta^-_{i+1,j+1} \\
		 & \zeta^-_{i-1,j}   & \zeta^-_{i,j}   & \zeta^-_{i+1,j}   \\
		 & \zeta^-_{i-1,j-1} & \zeta^-_{i,j-1} & \zeta^-_{i+1,j-1}
	\end{bmatrix} \\
	 & \mathbf{U}_{i,j} := \begin{bmatrix}
		 & u_{i-1,j+1} & u_{i,j+1} & u_{i+1, j+1} \\
		 & u_{i-1,j}   & u_{i,j}   & u_{i+1, j}   \\
		 & u_{i-1,j-1} & u_{i,j-1} & u_{i+1, j-1}
	\end{bmatrix},          & \mathbf{N}^\pm_{i,j} :=\begin{bmatrix}
		 & \omega_{i,j}^\pm(-1,1)  & \omega_{i,j}^\pm(0,1)  & \omega_{i,j}^\pm(1,1)  \\
		 & \omega_{i,j}^\pm(-1,0)  & 0                      & \omega_{i,j}^\pm(1,0)  \\
		 & \omega_{i,j}^\pm(-1,-1) & \omega_{i,j}^\pm(0,-1) & \omega_{i,j}^\pm(1,-1) \\
	\end{bmatrix}
\end{align*}
where $\mathbf{N^-}$ is a masking filter that passes the values in the negative neighborhood of node $(i,j)$.

We also introduce the Hadamard product $\odot$ between two identical matrices that creates another identical matrix with each entry being elementwise products. Moreover, double contraction of two tensors $A$ and $B$ is defined by $A : B = \sum A \odot B$ which is a scalar value and equals the sum of all entries of the Hadamard product of the tensors; \textit{i.e.}, note $A:A$ is square of Frobenius norm of $A$. Using these notations, the substitution rules read
\begin{equation}
	u_{i,j}^-=\begin{cases}
		u_{i,j}                                                                                                                                                                                                                                                                        & \text{ $\mathbf{x}_{i,j}\in \Omega^-$} \\
		\big(1 + \boldsymbol{\Gamma}^-_{i,j} : \mathbf{N}^-_{i,j}  \big) u_{i,j} -  \big( \boldsymbol{\Gamma}^-_{i,j}  \odot \mathbf{N}^-_{i,j} \big) : \mathbf{U}_{i,j}  - (\alpha + \delta_{i,j}\frac{\beta}{\mu^+}) \big(1 + \boldsymbol{\Gamma}^-_{i,j} : \mathbf{N}^-_{i,j} \big) & \text{ $\mathbf{x}_{i,j}\in \Omega^+$}
	\end{cases}
\end{equation}
\begin{equation}
	u_{i,j}^+=\begin{cases}
		\big(1 + \boldsymbol{\zeta}^-_{i,j} : \mathbf{N}^-_{i,j}   \big) u_{i,j} - \big( \boldsymbol{\zeta}^-_{i,j}  \odot \mathbf{N}^-_{i,j}  \big) : \mathbf{U}_{i,j} + \alpha + \delta_{i,j}\frac{\beta}{\mu^+} & \text{ $\mathbf{x}_{i,j}\in \Omega^-$} \\
		u_{i,j}                                                                                                                                                                                                    & \text{ $\mathbf{x}_{i,j}\in \Omega^+$}
	\end{cases}
\end{equation}

\begin{itemize}
	\item Rules based on approximating $\partial_\mathbf{n}u^-(\mathbf{r}^{pr}_{i,j})$:
\end{itemize}
\begin{equation}
	u_{i,j}^-=\begin{cases}
		u_{i,j}                                                                                                                        & \text{ $\mathbf{x}_{i,j}\in \Omega^-$} \\
		u_{i,j} (1-\zeta^+_{i,j}) - \sum_{(p,q)\in N_{i,j}^+} \zeta_{i,j,p,q}^+ u_{i+p,j+q} - \alpha - \delta_{i,j}\frac{\beta}{\mu^-} & \text{ $\mathbf{x}_{i,j}\in \Omega^+$}
	\end{cases}
\end{equation}
\begin{equation}
	u_{i,j}^+=\begin{cases}
		u_{i,j}(1 - \gamma^+_{i,j}) - \sum_{(p,q)\in N^+_{i,j}} \gamma^+_{i,j,p,q} u_{i+p, j+q} + (\alpha + \delta_{i,j}\frac{\beta}{\mu^-})(1 - \gamma^+_{i,j}) & \text{ $\mathbf{x}_{i,j}\in \Omega^-$} \\
		u_{i,j}                                                                                                                                                  & \text{ $\mathbf{x}_{i,j}\in \Omega^+$} \\
	\end{cases}
\end{equation}
in matrix notation we have
\begin{equation}
	u_{i,j}^-=\begin{cases}
		u_{i,j}                                                                                                                                                                                                 & \text{ $\mathbf{x}_{i,j}\in \Omega^-$} \\
		\big(1 + \boldsymbol{\zeta}^+_{i,j} : \mathbf{N}^+_{i,j}\big) u_{i,j}  -  \big( \boldsymbol{\zeta}^+_{i,j} \odot \mathbf{N}^+_{i,j} \big) : \mathbf{U}_{i,j} - \alpha - \delta_{i,j}\frac{\beta}{\mu^-} & \text{ $\mathbf{x}_{i,j}\in \Omega^+$}
	\end{cases}
\end{equation}
\begin{equation}
	u_{i,j}^+=\begin{cases}
		\big(1 + \boldsymbol{\Gamma}^+_{i,j} : \mathbf{N}^+_{i,j}\big) u_{i,j} - \big( \boldsymbol{ \Gamma}^+_{i,j} \odot  \mathbf{N}^+_{i,j}  \big) : \mathbf{U}_{i,j} + (\alpha + \delta_{i,j}\frac{\beta}{\mu^-})\big(1 + \boldsymbol{\Gamma}^+_{i,j} : \mathbf{N}^+_{i,j} \big) & \text{ $\mathbf{x}_{i,j}\in \Omega^-$} \\
		u_{i,j}                                                                                                                                                                                                                                                                     & \text{ $\mathbf{x}_{i,j}\in \Omega^+$} \\
	\end{cases}
\end{equation}
The overall algorithm is summarized in Algorithm \ref{euclid}.

\begin{algorithm}
	\caption{Bias Slow approximation of the non-existing solution value on a grid point based on existing solution values in its neighborhood. The notation is used for $u_{i,j}^\pm=B_{i,j}^\pm : \mathbf{U}_{i,j}+r_{i,j}^\pm$.}\label{euclid}
	\begin{algorithmic}[1]
		\Procedure{Bias Slow}{}
		\If {$\Gamma \cap \mathcal{C}_{i,j}= \emptyset$}
		\State $B_{i,j}^\pm=\begin{bmatrix}
				0 & 0 & 0 \\
				0 & 1 & 0 \\
				0 & 0 & 0 \\
			\end{bmatrix};\ \ r^\pm_{i,j}=0$
		\Else
		\If {$\mu_{i,j}^- > \mu_{i,j}^+$}
		\If {$\phi_{i,j}\ge 0$}
		\State $B_{i,j}^+=\begin{bmatrix}
				0 & 0 & 0 \\
				0 & 1 & 0 \\
				0 & 0 & 0 \\
			\end{bmatrix};\ \ r^+_{i,j}=0$
		\State $B_{i,j}^-=\begin{bmatrix}
				-\gamma_{i,j,-1,1}^-  & -\gamma_{i,j,0,1}^-  & -\gamma_{i,j,1,1}^-  \\
				-\gamma_{i,j,-1,0}^-  & 1-\gamma^-_{i,j}     & -\gamma_{i,j,1,0}^-  \\
				-\gamma_{i,j,-1,-1}^- & -\gamma_{i,j,0,-1}^- & -\gamma_{i,j,1,-1}^- \\
			\end{bmatrix};\ \ r^-_{i,j}=-(\alpha_{i,j}^{proj} + \delta_{i,j}\frac{ \beta_{i,j}^{proj}}{\mu_{proj}^+}) (1 - \gamma_{i,j}^-)$
		\Else

		\State $B_{i,j}^+=\begin{bmatrix}
				-\zeta_{i,j,-1,1}^-  & -\zeta_{i,j,0,1}^-  & -\zeta_{i,j,1,1}^-  \\
				-\zeta_{i,j,-1,0}^-  & 1-\zeta^-_{i,j}     & -\zeta_{i,j,1,0}^-  \\
				-\zeta_{i,j,-1,-1}^- & -\zeta_{i,j,0,-1}^- & -\zeta_{i,j,1,-1}^- \\
			\end{bmatrix};\ \ r^+_{i,j}=\alpha_{i,j}^{proj} + \delta_{i,j}\frac{ \beta_{i,j}^{proj}}{\mu_{proj}^+}$

		\State $B_{i,j}^-=\begin{bmatrix}
				0 & 0 & 0 \\
				0 & 1 & 0 \\
				0 & 0 & 0 \\
			\end{bmatrix};\ \ r^-_{i,j}=0$

		\EndIf

		\Else

		\If {$\phi_{i,j}\ge 0$}
		\State $B_{i,j}^+=\begin{bmatrix}
				0 & 0 & 0 \\
				0 & 1 & 0 \\
				0 & 0 & 0 \\
			\end{bmatrix};\ \ r^+_{i,j}=0$
		\State $B_{i,j}^-=\begin{bmatrix}
				-\zeta_{i,j,-1,1}^+  & -\zeta_{i,j,0,1}^+  & -\zeta_{i,j,1,1}^+  \\
				-\zeta_{i,j,-1,0}^+  & 1-\zeta^+_{i,j}     & -\zeta_{i,j,1,0}^+  \\
				-\zeta_{i,j,-1,-1}^+ & -\zeta_{i,j,0,-1}^+ & -\zeta_{i,j,1,-1}^+ \\
			\end{bmatrix};\ \ r^-_{i,j}=\alpha_{i,j}^{proj} + \delta_{i,j}\frac{ \beta_{i,j}^{proj}}{\mu_{proj}^-}$
		\Else

		\State  $B_{i,j}^+=\begin{bmatrix}
				-\gamma_{i,j,-1,1}^+  & -\gamma_{i,j,0,1}^+  & -\gamma_{i,j,1,1}^+  \\
				-\gamma_{i,j,-1,0}^+  & 1-\gamma^+_{i,j}     & -\gamma_{i,j,1,0}^+  \\
				-\gamma_{i,j,-1,-1}^+ & -\gamma_{i,j,0,-1}^+ & -\gamma_{i,j,1,-1}^+ \\
			\end{bmatrix};\ \ r^+_{i,j}=(\alpha_{i,j}^{proj} + \delta_{i,j}\frac{ \beta_{i,j}^{proj}}{\mu_{proj}^-}) (1 - \gamma_{i,j}^+)$

		\State $B_{i,j}^-=\begin{bmatrix}
				0 & 0 & 0 \\
				0 & 1 & 0 \\
				0 & 0 & 0 \\
			\end{bmatrix};\ \ r^-_{i,j}=0$

		\EndIf

		\EndIf

		\EndIf

		\EndProcedure
	\end{algorithmic}
\end{algorithm}

\subsection{Approach II. Finite discretization method fused with neural extrapolation}\label{sec:fusion}
We point out that although in approach I we used a regression-based method to impose the jump conditions on the grid points around the interface, it is possible to evaluate the neural network models as interpolation and extrapolation functions within the finite discretization scheme. Using the neural network models for solutions, we are able to evaluate extrapolations of the solution functions in a banded region around the interface as illustrated in figure \ref{fig:shapes}. Starting from the jump conditions, for points on the interface, $\mathbf{x}\in \Gamma$, we have
\begin{align*}
u^+ - u^- &= \alpha\\
\mu^+\partial_n u^+ - \mu^- \partial_n u^- &= \beta
\end{align*}
and after Taylor expansion in the normal direction we obtain on the adjacent grid points $(i,j)$
\begin{align}
u^+_{i,j} - u^-_{i,j}=[u]_{\mathbf{r}_{i,j}^{pr}} + \delta_{i,j}(\partial_\mathbf{n}u^+(\mathbf{r}^{pr}_{i,j}) - \partial_\mathbf{n}u^-(\mathbf{r}^{pr}_{i,j})) \label{eq::taylorexpandjump}
\end{align}
which explicitly incorporates the jump condition in the solutions. To incorporate the jump condition in fluxes we can rewrite either of the normal gradients in terms of the other 
\begin{align*}
\partial_\mathbf{n}u^+(\mathbf{r}^{pr}_{i,j}) &= \frac{\mu^-}{\mu^+}\partial_\mathbf{n}u^-(\mathbf{r}^{pr}_{i,j}) + \frac{\beta}{\mu^+}\\
\partial_\mathbf{n}u^-(\mathbf{r}^{pr}_{i,j}) &= \frac{\mu^+}{\mu^-}\partial_\mathbf{n}u^+(\mathbf{r}^{pr}_{i,j}) - \frac{\beta}{\mu^-} 
\end{align*}
which leads to two relationships among predictions of the two neural networks at each grid point in the banded extrapolation region
\begin{align}
u^+_{i,j} - u^-_{i,j}=\alpha(\mathbf{r}_{i,j}^{pr}) + \delta_{i,j}\bigg(\big(\frac{\mu^-}{\mu^+}-1\big)\partial_\mathbf{n}u^-(\mathbf{r}^{pr}_{i,j}) + \frac{\beta(\mathbf{r}_{i,j}^{pr})}{\mu^+} \bigg)  \label{eq::extrapolate1}\\
u^+_{i,j} - u^-_{i,j}=\alpha(\mathbf{r}_{i,j}^{pr}) + \delta_{i,j}\bigg(\big(1-\frac{\mu^+}{\mu^-}\big)\partial_\mathbf{n}u^+(\mathbf{r}^{pr}_{i,j}) + \frac{\beta(\mathbf{r}_{i,j}^{pr})}{\mu^-} \bigg) \label{eq::extrapolate2}
\end{align}

Note that we are representing solution functions, $\hat{u}^\pm(\mathbf{r})$, with neural networks where computing the normal derivatives is trivial using automatic differentiation of the network along the normal directions. In contrast to finite discretization methods, solutions at off-grid points is readily available by simply evaluating the neural network function at any desired points. Note that we can compute the projected location on the interface starting from each grid point $(i,j)$ using the level-set function:
\begin{align*}
\mathbf{r}_{ij}^{proj}=\mathbf{r}_{ij} - \delta_{i,j} \mathbf{n}_{i,j}
\end{align*}

In the second approach, the loss function remains as before, except the unknown $u^\pm$ values are derived using equations \ref{eq::extrapolate1}--\ref{eq::extrapolate2}, instead of computing a regression-based extrapolation function based on the points in the neighborhood of interface cells: 

\begin{align*}
\mathcal{L} = & \bigg\vert\bigg\vert \sum_{s=-,+} k_{i,j}^s u_{i,j}^{s} |\mathcal{V}_{i,j}^s| - \sum_{s=-,+}\bigg( \mu_{i-\frac{1}{2},j}^s A_{i-\frac{1}{2},j}^s\frac{u_{i-1,j}^s - u_{i,j}^s}{\Delta x}     +   \mu_{i+\frac{1}{2},j}^s A_{i+\frac{1}{2},j}^s\frac{u_{i+1,j}^s - u_{i,j}^s}{\Delta x} +                                                          \\
	 & \mu_{i, j-\frac{1}{2}}^s A_{i, j-\frac{1}{2}}^s\frac{u_{i,j-1}^s - u_{i,j}^s}{\Delta y} + \mu_{i, j+\frac{1}{2}}^s A_{i, j+\frac{1}{2}}^s\frac{u_{i,j+1}^s - u_{i,j}^s}{\Delta y} \bigg)  -  \sum_{s=-,+} f_{i,j}^{s} |\mathcal{V}_{i,j}^s| - \int_{\Gamma\cap \mathcal{V}_{i,j}} \beta d\Gamma \bigg\vert\bigg\vert_2^2
\end{align*}

However, there is a major downside with this approach for training because the automatic differentiation has to be applied on the network once more that effectively amounts to compute second-order derivatives of the network. This slows down convergence, and the time-to-solution increases with square of depth of the neural network while in the regression-based method the cost grows linearly in the network depth by restricting to only first order automatic differentiation.

\subsection{Optimization scheme}

One of the main benefits of NBM is enabling the application of techniques from the vast literature on preconditioning linear systems to accelerate training of neural network models for the solution of PDEs. We note that in NBM these preconditioners do not need to be differentiable as long as their operations only depend on the geometry and physical properties of the domains, and not explicitly on the solution values of the PDE. Therefore existing software libraries for preconditioning could be used in \texttt{JAX-DIPS}. In this section we introduce the optimization techniques for training neural network models in \texttt{JAX-DIPS}.

\subsubsection{Preconditioners are ideal network regularizers for solving PDEs}
Finite discretization methods lead to solving a linear algebraic system with gaurantees on convergence and accuracy. The geometric irregularities and fine-grain details of the system around interfaces often lead to bad condition number for the linear system, which can be remedied by applying preconditioners. Intuitively, condition number is caused by a separation of scales for geometric lengthscales or material properties that underly the solution patterns. One of the strengths of the presented approach is to readily enable usage of preconditioners for training neural network surrogate models.

Preconditioners are a powerful technique to accelerate convergence of tradional numerical linear algebraic solvers. Given a poorly conditioned linear system $Ax=b$ one can obtain an equivalent system $\hat{A}\hat{x}=\hat{b}$ with accelerated convergence rate when using iterative gradient based methods. For the conjugate gradient method convergence iteration is proportional to $\sqrt{\kappa(A)}$ where $\kappa(A)$ is the condition number of matrix $A$. Preconditioning is achieved by mapping the linear system with a nonsingular matrix $M$ into a new space $M^{-1} A x = M^{-1}b$ where $M^{-1}A$ has more regular spread of eigenvalues, hence a better condition number. The precondition matrix $M$ should approximate $A^{-1}$ such that $\vert I - M^{-1}A\vert < 1$. The simplest choice is the Jacobi prconditioner which amounts to using the diagonal part of $A$ as the preconditioner, $M=diag(A)$. Note that the diagonal term is locally available at each point and it is straightforward to parallelize.

In this work we use the Jacobi pre-conditioner. Basically, every element of the left-hand-side ($A u$) and right-hand-side ($b$) vectors are divided by the coefficient of the diagonal term of the matrix given by:
\begin{align*}
a_{ii} = \sum_{s=-,+} \bigg( k_{i,i}^s \vert \mathcal{V}_{i,i}^s\vert + (\mu_{i-\frac{1}{2},i}^s A_{i-\frac{1}{2},i}^s +  \mu_{i+\frac{1}{2},i}^s A_{i+\frac{1}{2},i}^s ) / \Delta x 
+  (\mu_{i, i-\frac{1}{2}}^s A_{i, i-\frac{1}{2}}^s +  \mu_{i, i+\frac{1}{2}}^s A_{i, i+\frac{1}{2}}^s ) / \Delta y \bigg)
\end{align*}
Note that for memory efficiency we never explicitly compute the matrix, instead we compute the effect of matrix product of $A u$.

\subsubsection{Learning rate scheduling}
First order methods have longer time-to-solution but require less memory, while second order methods are faster to converge but require massive memory footprint. In \texttt{JAX-DIPS} we primarily utilize first order optimization methods such as Adam \cite{kingma2014adam}. Second order methods such as Newton or BFGS certainly offer convergence in less iterations but require much more memory. Traditionally used GMRES or Conjugate Gradient methods for sparse linear systems are somewhere between first order and second order optimization methods that are based on building basis vectors by computing gradients that are conjugate to each other $\mathbf{p}_j^T A \mathbf{p}_i=0$ and will converge to the solution in at most $n$ steps; \textit{i.e.}, at most the solution vector is spanned in the full basis. We found that starting from a zero guess for the solution it is important to start from a large learning rate and gradually decay the learning rate in a process of exponential annealing. For this purpose, we use the exponential decay scheduler provided by \texttt{Optax} \cite{optax2020github} to control the learning rate in the Adam optimizer:
\begin{align*}
r_{k} = r_0 \alpha^{k / T}
\end{align*}
where $r_k$ is the learning rate at step $k$ of optimization, $\alpha<1$ is the decay rate, and $T$ is the decay count-scale. By default, we set $T=100$, $\alpha=0.975$, starting from an initial value of $r_0=10^{-2}$ and clip gradients by maximum global gradient norm (to a value $1$) \cite{gradClipping} before applying the Adam updates in each step. We note a larger decay rate, \textit{e.g.} $\alpha=0.98$, leads to small oscillations after $10000$ steps and although similar levels of accuracy can be achieved at much less iterations, here we report results with the more robust decay rate. 


\subsubsection{Domain switching optimization scheme}
The linear system suffers from worse condition number in the domain with more variability in diffusion coefficient, or where diffusion coefficient is larger; \textit{i.e.,} the fast region. This leads to regionally unbalanced solution error where the overall error is systematically lopsided by the faster diffusion region. We found this problem can be improved by interleaving region-specific optimization epochs in the training pipeline, where only one of the networks is updated based on the loss computed in its region. See Algorithm \ref{switching} for details of the algorithm. We note that an alternative strategy for addressing this issue is scaling the gradients for the fast region with respect to the slow region during training.

\begin{algorithm}
	\caption{Domain switching method. Switching interval is $\tau$.}\label{switching}
	\begin{algorithmic}[1]
		\Procedure{Domain Switching Optimization}{}
		
		\For {$epoch$ in $0\cdots$ N}	
			\State $region = Region (epoch)$
			
			\If {$region > 0$}
				\If {$\mu^- > \mu^+$}
					\State optimize $u^-_{NN}$ in $\Omega^-$ given fixed $u^+_{NN}$  
				\Else
					\State optimize $u^+_{NN}$ in $\Omega^+$ given fixed $u^-_{NN}$ 
				\EndIf	
			\EndIf	
			\\
			\If {$region==0$}
				\State optimize both networks in $\Omega^- \cup \Omega^+$ 
			\EndIf	
			\\
			\If {$region < 0$}
				\If {$\mu^- < \mu^+$}
					\State optimize $u^-_{NN}$ in $\Omega^-$ given fixed $u^+_{NN}$  
				\Else
					\State optimize $u^+_{NN}$ in $\Omega^+$ given fixed $u^-_{NN}$ 
				\EndIf					
			\EndIf

		\EndFor

		\EndProcedure
		\\
		\Procedure{Region}{epoch}
			\If {$mode == \textrm{whole region} \rightarrow \textrm{fast region}$}
				\State $region = epoch\ \%\ \tau $ 
			\EndIf
			\If {$mode == \textrm{fast region}\rightarrow \textrm{whole region} \rightarrow \textrm{slow region}$}
				\State $region = \tau // 2 - epoch\ \%\ \tau  $ 
			\EndIf
		\EndProcedure
	\end{algorithmic}
\end{algorithm}

\subsubsection{Multi-GPU parallelization with data parallel training}
The NBM is embarrasingly parallel and residual evaluation at each point is independent from other points. Therefore, multi-GPU parallelization does not involve inter-GPU communication for evaluating the residuals per point. We partition the training points and distribute them along with copies of neural network parameters among multiple GPUs to compute gradient updates per batch. Then we aggregate these updates by averaging the values on different GPUs. The updates are then broadcasted and model parameters are updated on each device \cite{shallue2018measuring}.

The ability to batch over grid points is one of the key enabling factors for reaching higher resolutions and higher dimensions, as it allows to set a limit on the required GPU memory. With NBM it is straightforward to scale finite discretization methods on GPU clusters. 

\section{Numerical Results} \label{sec:examples}
We consider examples for solution to elliptic problems of the form
\begin{align*}
	 & k^{\pm}u^{\pm} - \nabla \cdot (\mu^{\pm}\nabla u^\pm)=f^{\pm}, & \mathbf{x}\in\Omega^\pm \\
	 & [u]=\alpha,                                                    & \mathbf{x} \in \Gamma   \\
	 & [\mu \partial_{\mathbf{n}}u]=\beta,                            & \mathbf{x} \in \Gamma
\end{align*}

Using different features of \texttt{JAX-DIPS} one can compose solvers with different training configurations; \textit{i.e.}, single/multi-resolution, single/multi-batch, and single/multi-GPU, and domain alternating training. Moreover, the neural extrapolation method discussed in section \ref{sec:fusion} provides an alternative solver. Below we implement and compare numerical accuracy and performance of these strategies. 

For each accuracy metric we report \textit{order of convergence}. Order of convergence, denoted by $\rm p$, is computed by doubling the number of grid points in every dimension and measuring the $\rm L^\infty$ error of solution and its gradient over all the grid points in the domain:
\begin{align*}
&\frac{\textrm{err}(2h)}{\textrm{err}(h)}=2^p \rightarrow p = \log_2\bigg(\frac{\textrm{err}(2h)}{\textrm{err}(h)}\bigg)  & h=\min(h_x,h_y,h_z)
\end{align*}

\subsection{Accuracy in the bulk: no interface}\label{sec:ex1}
As baseline we consider the solution in the bulk at the absence of interfaces. Computational domain is $\Omega\in [-1,1]^3$ with the exact solution given by $u(x,y,z)=\cos(x)\sin(y)\cos(z)$, coefficients $\mu=1$ and $k=0$ and the source term is $f(x,y,z)=3\cos(x)\sin(y)\cos(z)$. The accuracy results are reported in Table \ref{tab:convergence_bulk}. The neural network consists of 5 hidden layers with 10 neurons with \texttt{CeLU} activation function, consisting of 491 total trainable parameters. In each case the batch size is equal to the number of training grid points in order to ensure only a single batch training. We note this setting identifies the maximum memory usage and the minimum time per epoch. Decreasing the batch size reduces the amount of GPU memory that is required by folding the training data into several passes of the optimizer at the expense of increasing the time needed per epoch. This degree of freedom accomodates for adjusting to the available hardware specifications. 

To quantify the solution error we construct a high resolution \textit{evaluation grid} with $256\times 256 \times 256$ grid points to account for spatial regions outside of the training grid. After training the network with specified training resolution we evaluate the network over the evaluation grid and compute $\rm RMSE$, $\rm L^\infty$, and relative $\rm L^2$ errors. We also report the GPU utilization during our experiments. It is also important to note the number of epochs determines the level of accuracy in our experiments, here we wait $10,000$ epochs before measuring the errors.
\begin{table}[ht]
\begin{center}
\caption{Convergence test of example \ref{sec:ex1}. Timings include the initial compilation time. Measurements are on a single NVIDIA A6000 GPU. The model has 5 hidden layers with 10 neurons each, overall 491 trainable parameters.} \label{tab:convergence_bulk}
\begin{adjustbox}{width=\textwidth}
\begin{tabular}{ccccccccc}
\toprule
 & \multicolumn{2}{c}{RMSE}& \multicolumn{2}{c}{$\rm L^\infty$} & \multicolumn{2}{c}{Rel. $\rm L^2$} & \multicolumn{2}{c}{GPU Statistics} \\
\cmidrule(r){2-3}\cmidrule(r){4-5}\cmidrule(r){6-7} \cmidrule(r){8-9}
$\rm N_{x,y,z}$   &   Solution    &   Order   &   Solution   &   Order &  Solution   &   Order & t (sec/epoch) & VRAM (GB)\\
\midrule
$2^3$ & $1.07\times 10^{-1}$ &    $-$  & $2.13\times 10^{-2}$  &    $-$   & $2.42\times 10^{-2}$ &   $-$   & $0.012$ & $0.81$ \\
$2^4$ & $8.04\times 10^{-3}$ &  $3.73$ & $1.67\times 10^{-2}$  &  $0.35$  & $1.80\times 10^{-2}$ & $0.43$  & $0.024$ & $2.22$ \\ 
$2^5$ & $1.83\times 10^{-3}$ &  $2.13$ & $5.24\times 10^{-3}$  &  $1.67$  & $4.09\times 10^{-3}$ & $2.13$  & $0.091$ & $5.77$ \\ \bottomrule
\end{tabular}
\end{adjustbox}
\end{center}
\end{table}

\subsection{Accuracy on spherical interface: single-resolution, single batch, single GPU}
We use a single uniform grid and train on all the points in a single batch. We consider a sphere $\phi(\mathbf{x})=\sqrt{x^2 + y^2 + z^2} - 0.5$ centered in a box $\Omega:[-1,1]^3$ with the exact solution
\begin{align*}
& u^-(x,y,z)=e^{z}, & \phi(\mathbf{x})<0\\
& u^+(x,y,z)=\cos(x)\sin(y), & \phi(\mathbf{x})\ge 0
\end{align*}
and variable diffusion coefficients
\begin{align*}
\mu^-(x,y,z)&=y^2 \ln(x+2) + 4 &\phi(\mathbf{x})<0 \\
\mu^+(x,y,z)&=e^{-z} &\phi(\mathbf{x})\ge 0 
\end{align*}
that imply variable source terms
\begin{align*}
&f^-(x,y,z)=-[y^2\ln(x+2) + 4] e^{z} &\phi(\mathbf{x})< 0\\
&f^+(x,y,z)=2\cos(x)\sin(y)e^{-z} &\phi(\mathbf{x})\ge 0
\end{align*}

The network has 5 hidden layers with 10 neurons in each layer using sine activation functions. Table \ref{tab:convergence_sphere} reports convergence results for the solution in the $L^\infty$-norm and root-mean-squared-error of the solution. 

\begin{figure}[ht]
	\centering
	\includegraphics[width=0.49\linewidth]{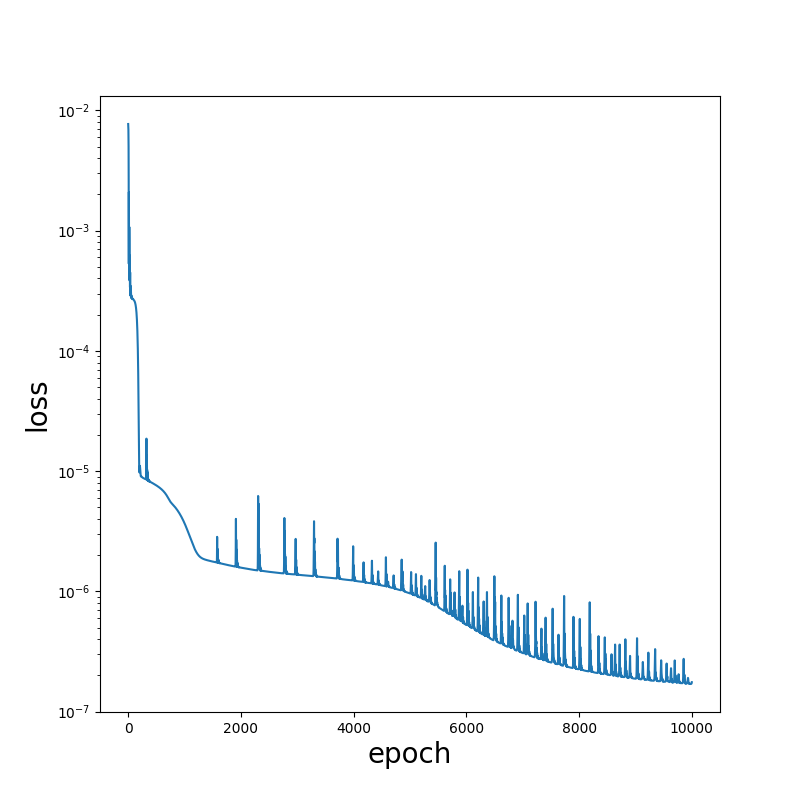}
	\includegraphics[width=0.45\linewidth]{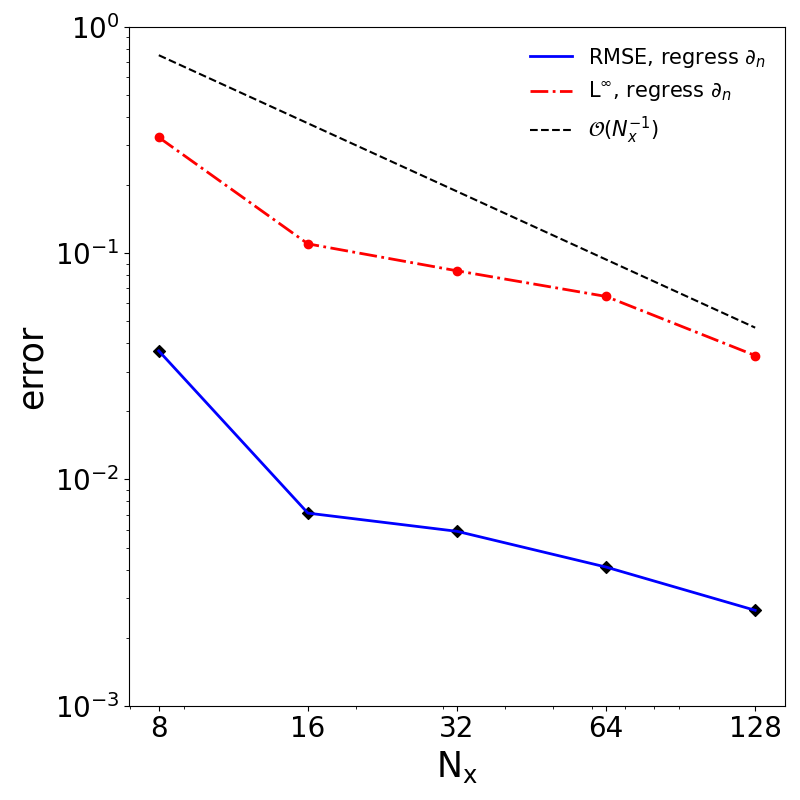}
	\caption{Loss evolution with epochs for the sphere of $16\times 16\times 16$ grid (left) and different accuracy measures, RMSE and $L^\infty$, at 5 different resolutions (right).}
	\label{fig:losses}
\end{figure}

\begin{table}[ht]
\begin{center}
\caption{Convergence and timings for the sphere example averaged over $10,000$ epochs. Timings include the initial compilation time. Measurements are on a single NVIDIA A6000 GPU. The regression-based method has 5 hidden layers with 10 neurons each, overall 982 trainable parameters.} \label{tab:convergence_sphere}
\begin{adjustbox}{width=\textwidth}
\begin{tabular}{ccccccc}
\toprule
 & \multicolumn{2}{c}{RMSE}& \multicolumn{2}{c}{$L^\infty$} & \multicolumn{2}{c}{GPU Statistics} \\
\cmidrule(r){2-3}\cmidrule(r){4-5}\cmidrule(r){6-7}
$\rm N_{x,y,z}$   &   Solution    &   Order   &   Solution   &   Order & t (sec/epoch) & VRAM (GB)\\
\midrule
$2^3$ & $3.7\times 10^{-2}$ &  -        & $3.25\times 10^{-1}$  &   -     & $0.0306$ & $1.05$ \\
$2^4$ & $7.1\times 10^{-3}$ &  $2.38$   & $1.10\times 10^{-1}$   & $1.56$  & $0.056$  & $1.72$ \\ 
$2^5$ & $5.9\times 10^{-3}$ &  $0.27$   & $8.36 \times 10^{-2}$ & $0.4$   & $0.053$  & $2.15$ \\ 
$2^6$ & $4.1\times 10^{-3}$ &  $0.53$   & $6.44\times 10^{-2}$  & $0.38$  & $0.287$  & $5.57$ \\ 
$2^7$ & $2.64\times 10^{-3}$&  $0.64$   & $3.53\times 10^{-2}$  & $0.87$  & $2.125$ & $32.1$ \\ \bottomrule
\end{tabular}
\end{adjustbox}
\end{center}
\end{table}

%

\begin{figure}
     \centering
     \begin{subfigure}[b]{\textwidth}
         \centering
         \includegraphics[width=\linewidth]{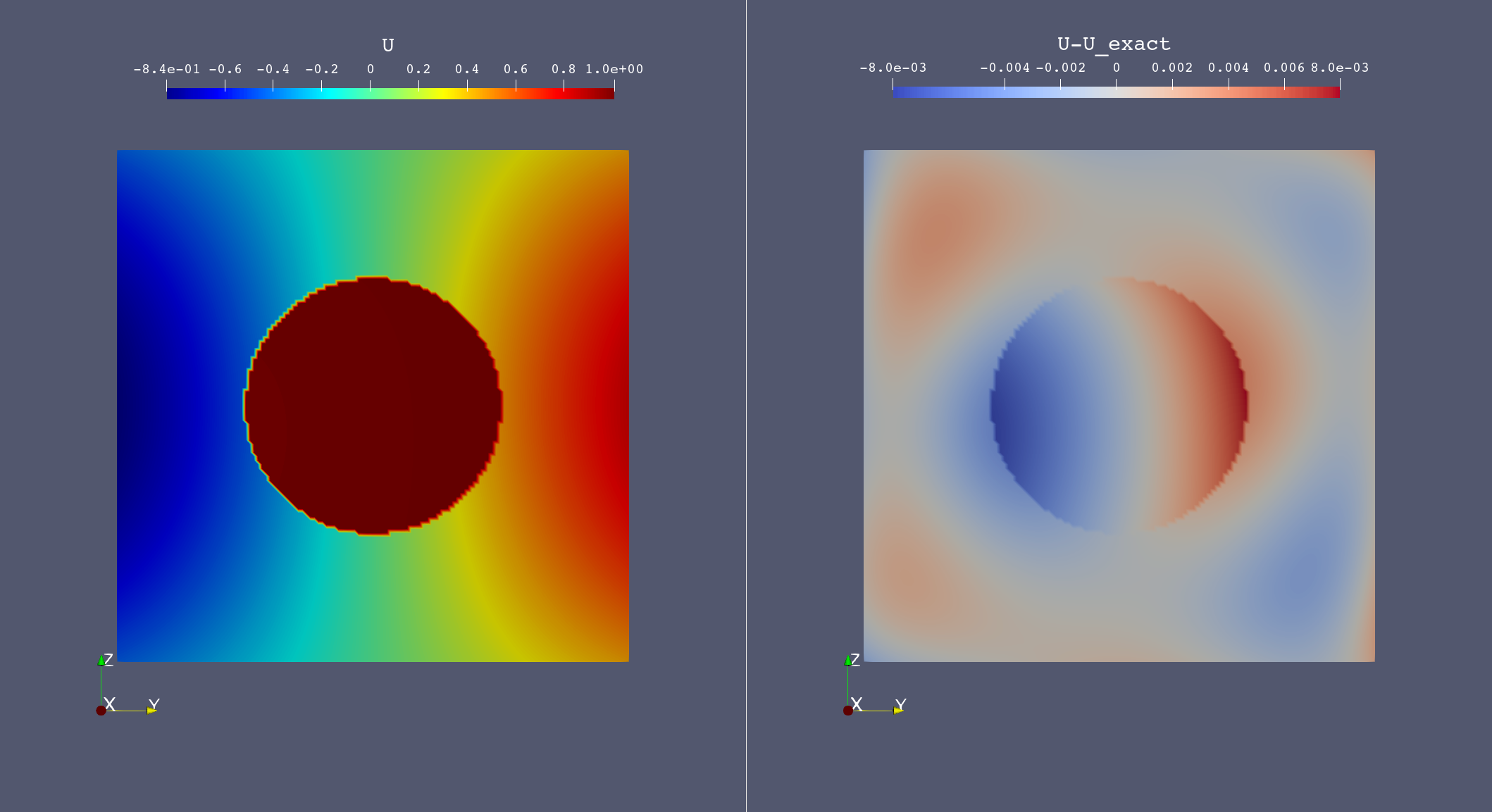}
	\caption{Illustration of numerical solution and absolute error on a cross section of the domain. }
	\label{subfig:sphere}
     \end{subfigure}
     \begin{subfigure}[b]{\textwidth}
         \centering
         \includegraphics[width=\linewidth]{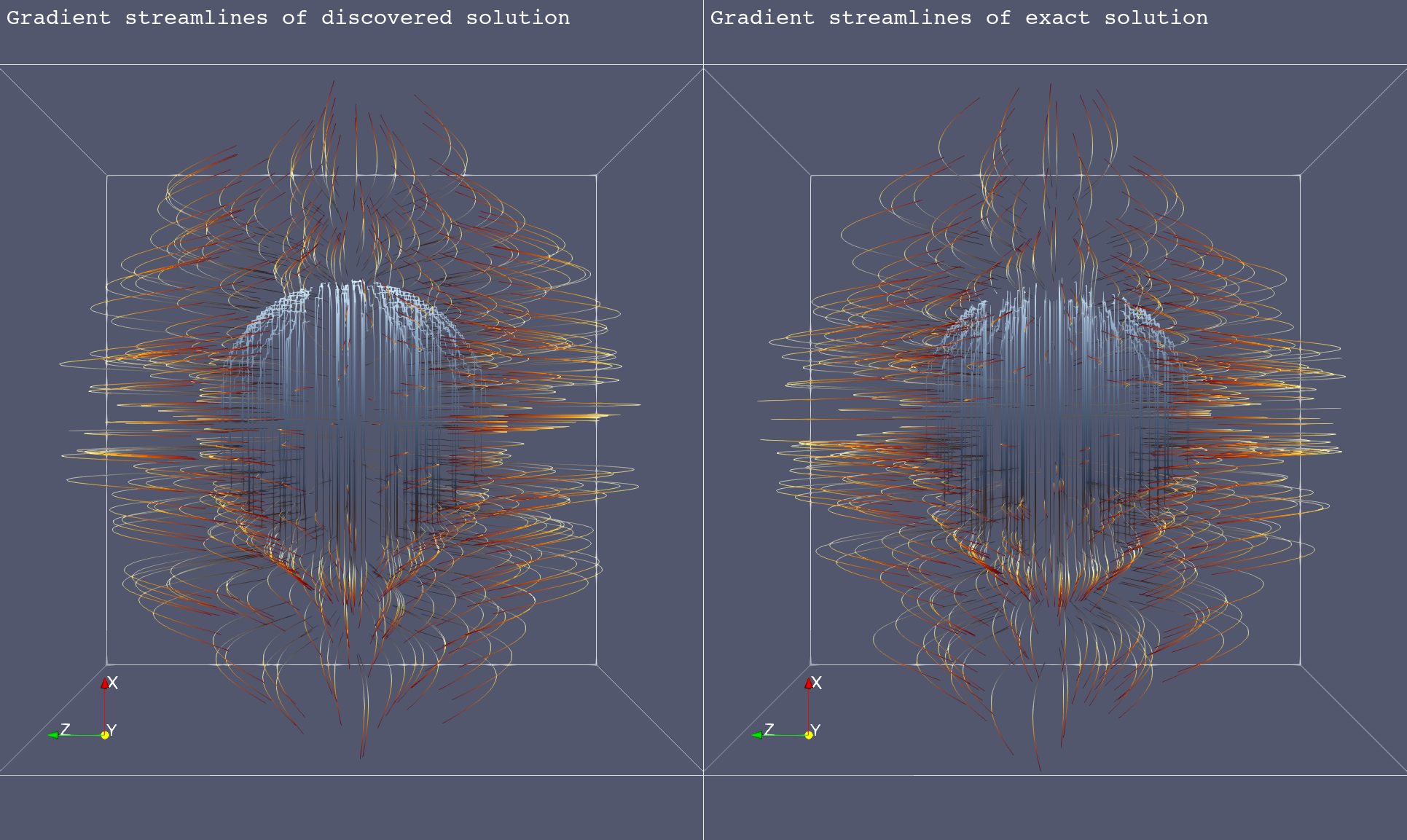}
	\caption{Streamlines of solution gradient for (left) the surrogate neural model colored by model solution value, (right) exact streamlines colored by exact solution values.  }
         \label{fig:spheregrad}
     \end{subfigure}
        \caption{The neural network surrogate model trained on a $128^3$ grid using a single NVIDIA A6000 GPU.}
        \label{fig:sphere}
\end{figure}



\subsection{Accuracy on star interface: single GPU, domain switching, neural extrapolation, and batching} \label{subsec:star}

We consider a more irregular interface in this section where the less favorable condition number of the discretized linear system is known to degrade the accuracy of results in traditional solvers. Figures 5 and 6 in \cite{bochkov2020solving} characterize degradation of solution errors due to worse condition number in the method bootstrapped in our work. We note this is a main challenge facing all finite discretization schemes for interfacial PDE problems.

For the star example we use a pair of fully connected feedforward neural networks, each composed of $1$ hidden layer and $100$ neurons with \texttt{sine} activation function, followed  by an output layer with $1$ linear neuron. There are a total of $1,002$ trainable parameters in the model. We consider a star-shaped interface with inner and outer radii $r_i=0.151$ and $r_e=0.911$ that is immersed in a box $\Omega:[-1,1]^3$ described by the level-set function
\begin{align*}
\phi(\mathbf{x}) = \sqrt{x^2 + y^2 + z^2} - r_0 \bigg( 1 + \big( \frac{x^2 + y^2}{x^2 + y^2 + z^2}\big)^2  \sum_{k=1}^3 \beta_k \cos\big(n_k \big(\arctan\big(\frac{y}{x}\big) - \theta_k\big) \big) \bigg)
\end{align*}
with the parameters
\begin{align*}
&r_0 = 0.483, &\begin{pmatrix}
n_1\\
\beta_1\\
\theta_1
\end{pmatrix}=\begin{pmatrix}
3\\
0.1\\
0.5
\end{pmatrix}, &\begin{pmatrix}
n_2\\
\beta_2\\
\theta_2
\end{pmatrix}=\begin{pmatrix}
4\\
-0.1\\
1.8
\end{pmatrix}, &\begin{pmatrix}
n_3\\
\beta_3\\
\theta_3
\end{pmatrix}=\begin{pmatrix}
7\\
0.15\\
0
\end{pmatrix}
\end{align*}
Considering an exact solution 
\begin{align*}
& u^-(x,y,z)=\sin(2x)\cos(2y) e^{z}, & \phi(\mathbf{x})<0\\
& u^+(x,y,z)=\bigg[ 16\big(\frac{y-x}{3}\big)^5 - 20 \big(\frac{y-x}{3}\big)^3 + 5\big( \frac{y-x}{3}\big) \bigg]\ln(x+y+3)\cos(z)          , & \phi(\mathbf{x})\ge 0
\end{align*}
and the diffusion coefficient
\begin{align*}
\mu^-(x,y,z)&=10\bigg[ 1+0.2\cos(2\pi(x+y))\sin(2\pi(x-y)) \cos(z) \bigg] &\phi(\mathbf{x})<0 \\
\mu^+(x,y,z)&=1 &\phi(\mathbf{x})\ge 0 
\end{align*}
\begin{table}[ht]
\begin{center}
\caption{Convergence in solution of the star geometry using the single-resolution regression-based solver with domain switching. We report $L^\infty$-norm error as well as root-mean-squared-error (RMSE) of the solution field evaluated everywhere in the domain. Timings are averaged over $10,000$ epochs in each case and include the initial compilation time for jaxpressions. The neural network pair have 1 hidden layer each with 100 neurons, overall $1,002$ trainable parameteres. Domain switching scheme follows the $\textrm{whole region} \rightarrow \textrm{fast region} \rightarrow \textrm{fast region}$ sequence.}
\label{tab:convergence_star}
\begin{adjustbox}{width=\textwidth}
\begin{tabular}{ccccccc}
\toprule
 & \multicolumn{2}{c}{RMSE}& \multicolumn{2}{c}{$L^\infty$} & \multicolumn{2}{c}{GPU Statistics} \\
\cmidrule(r){2-3}\cmidrule(r){4-5}\cmidrule(r){6-7} 
$\rm N_{x,y,z}$   &   Solution    &   Order   &   Solution   &   Order & t (sec/epoch) & VRAM (GB)\\
\midrule
\rowcolor{lightgray}\multicolumn{7}{c}{regress $\partial_n$}  \\
$2^3$ & $ 1.36\times 10^{-1}$ &  -         & $ 1.27$                 &   -   & $ 0.019  $  & $  0.98 $ \\
$2^4$ & $ 7.98\times 10^{-2}$ &  $0.77 $   & $ 8.23\times 10^{-1} $  & $ 0.63$   & $ 0.022  $  & $ 1.01  $ \\ 
$2^5$ & $ 4.36\times 10^{-2}$ &  $0.87 $   & $ 3.85\times 10^{-1} $  & $ 1.10 $  & $ 0.032  $  & $ 1.30  $ \\ 
$2^6$ & $ 2.43\times 10^{-2}$ &  $0.84 $   & $ 2.28\times 10^{-1}$   & $ 0.76$   & $ 0.200  $  & $ 3.7 $ \\  
\rowcolor{lightgray}\multicolumn{7}{c}{neural $\partial_n$}  \\
$2^3$ & $2.17\times 10^{-1} $ &      -     & $2.89 $                &   -      & $0.0259 $ & $ 0.93$ \\
$2^4$ & $1.34\times 10^{-1} $ &  $0.70 $   & $ 1.66$                & $ 0.80 $ & $ 0.0408$ & $1.19 $ \\ 
$2^5$ & $5.68\times 10^{-2} $ &  $1.24 $   & $ 8.17\times 10^{-1}$  & $ 1.02$  & $ 0.0712$ & $ 2.96$ \\ 
$2^6$ & $ 2.77\times 10^{-2}$ &  $1.03 $   & $3.94\times 10^{-1} $  & $ 1.05$   & $ 0.334 $ & $ 13.6$ \\  \bottomrule
\end{tabular}
\end{adjustbox}
\caption{Convergence in solution of the star geometry using the multi-resolution regression-based solver with batching. 
	We use a multi-resolution training protocol that refines to 4 levels at each collocation point, this slightly improves accuracies although in the current version of \texttt{JAX-DIPS} ($\rm v0.0.1$) the memory requirement increases. Batch size is the minimum of $64\times 64\times 32$ and number of collocation points, which ensures memory saturation at $~30$ GB.}
	\label{tab:convergence_star_batching}
	\begin{adjustbox}{width=\textwidth}
	\begin{tabular}{ccccccc}
	\toprule
	regress $\partial_n$ & \multicolumn{2}{c}{RMSE}& \multicolumn{2}{c}{$L^\infty$} & \multicolumn{2}{c}{GPU Statistics} \\
	\cmidrule(r){1-1}\cmidrule(r){2-3}\cmidrule(r){4-5}\cmidrule(r){6-6}\cmidrule(r){7-7}	
	$\rm N_{x,y,z}$   &   Solution    &   Order   &   Solution   &   Order & t (sec/epoch) & VRAM (GB)\\
	\midrule
	$2^3$ & $ 1.05\times 10^{-1}$ &   -        & $ 1.29 $                &   -       & $ 0.0225 $     & $ 1.27 $ \\
	$2^4$ & $ 5.52\times 10^{-2}$ &  $0.93 $   & $ 6.22\times 10^{-1} $  & $ 1.05$   & $ 0.0411 $     & $ 1.27 $ \\ 
	$2^5$ & $ 2.44\times 10^{-2}$ &  $1.18 $   & $ 2.66\times 10^{-1}$   & $ 1.23$   & $ 0.1814 $     & $ 8.3  $ \\ 
	$2^6$ & $ 2.33\times 10^{-2}$ &  $0.07 $   & $ 2.24\times 10^{-1} $  & $ 0.25$   & $ 1.889 $      & $ 29.6 $ \\ 
	\rowcolor{pink}$2^7$ & $ 8.62\times 10^{-2}$ &  $-1.88 $   & $ 3.80\times 10^{-1} $   & $-0.76$   & $ 9.649 $  & $ 29.7 $ \\ \bottomrule
	\end{tabular}
	\end{adjustbox}
\end{center}
\end{table}
\begin{figure}
     \centering
     \begin{subfigure}[b]{\textwidth}
         \centering
         \includegraphics[width=0.49\linewidth]{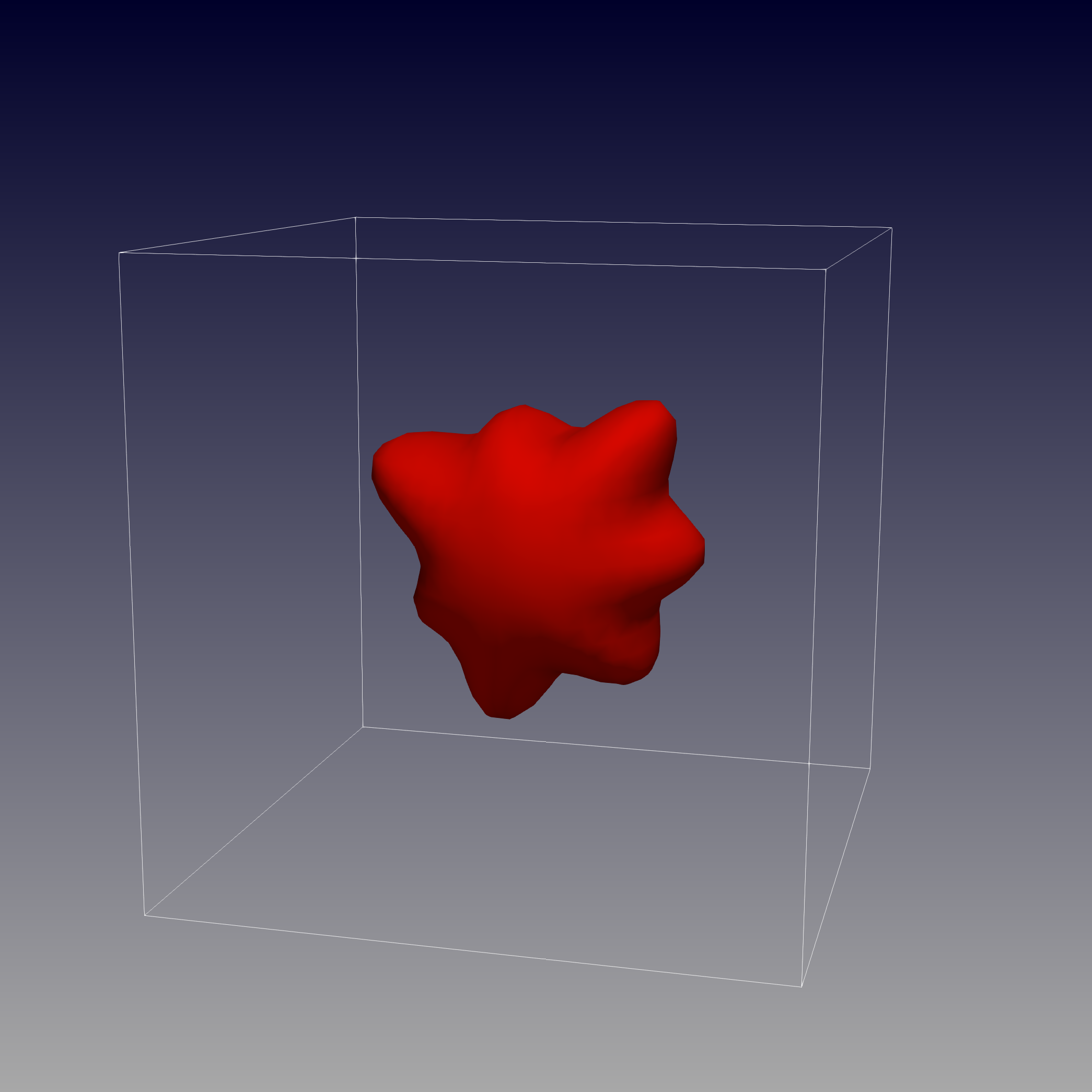}
		 \includegraphics[width=0.49\linewidth]{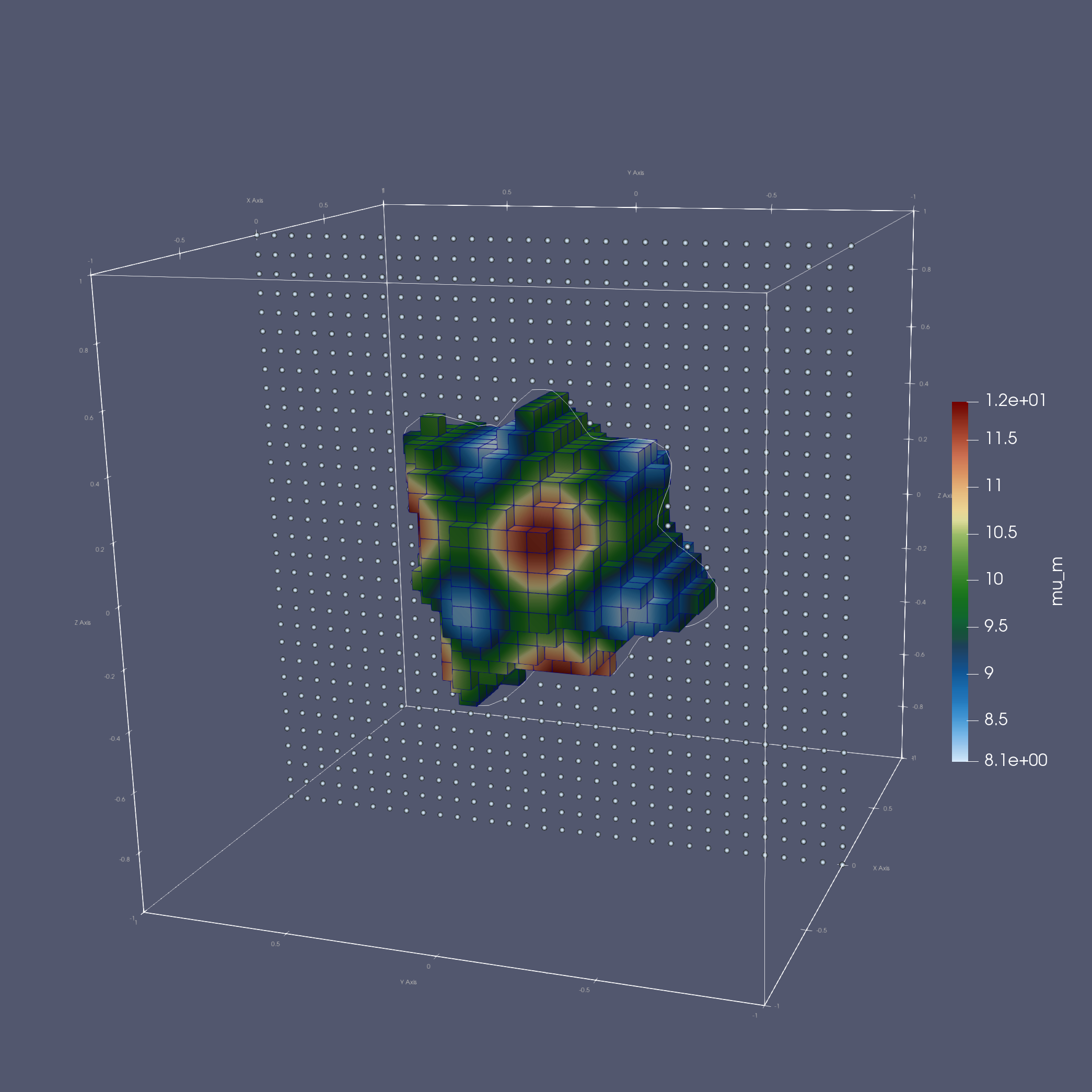}
		 \caption{Illustration of three dimensional interface used (left), and $\mu^\pm$ on the $32\times 32\times 32$ grid (right).}
		 \label{subfig:star}
     \end{subfigure}
     \begin{subfigure}[b]{\textwidth}
         \centering
         \includegraphics[width=0.49\linewidth]{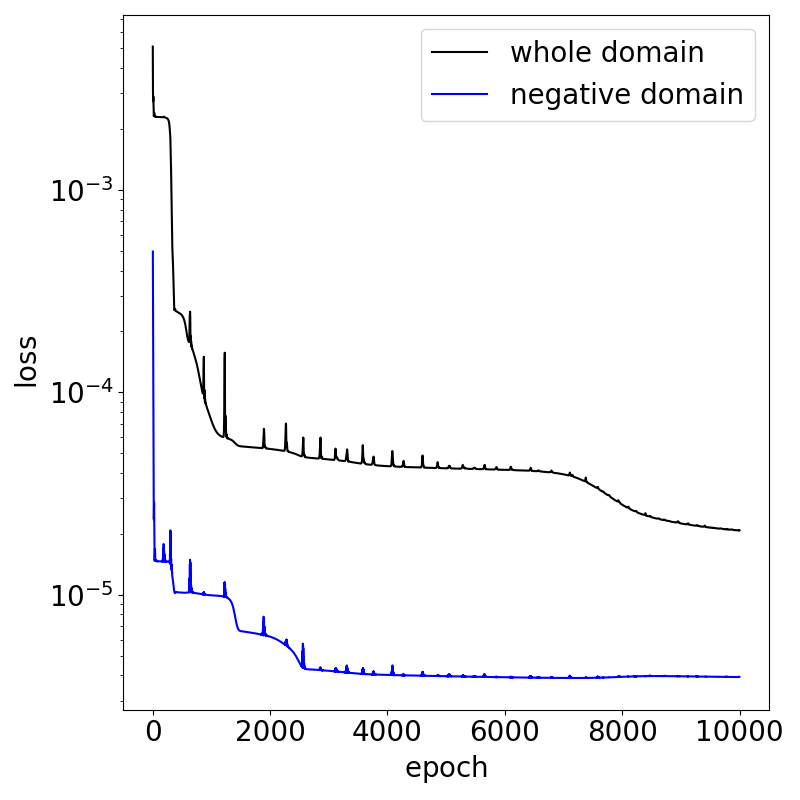}
		 \includegraphics[width=0.49\linewidth]{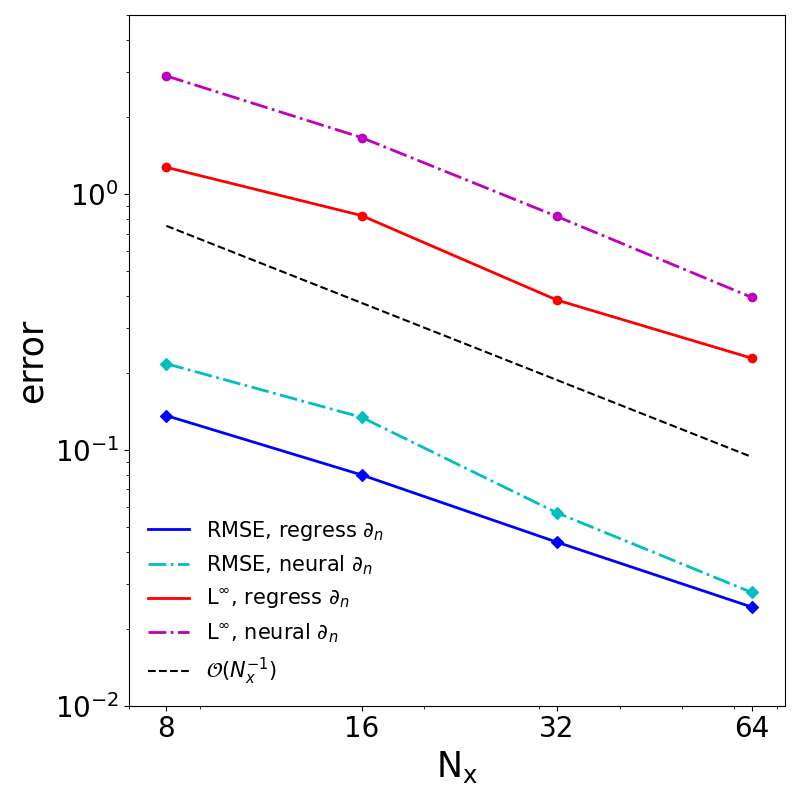}
		 \caption{Loss evolution with epochs for the star of $64\times 64\times 64$ grid using domain switching training (left), and decrease in error by increasing resolutions (right).}
		 \label{fig:lossestar}
     \end{subfigure}
        \caption{The neural network model trained with different configurations and resolutions.}
        \label{fig:star}
\end{figure}
\begin{figure}
	\centering
	\includegraphics[width=\linewidth]{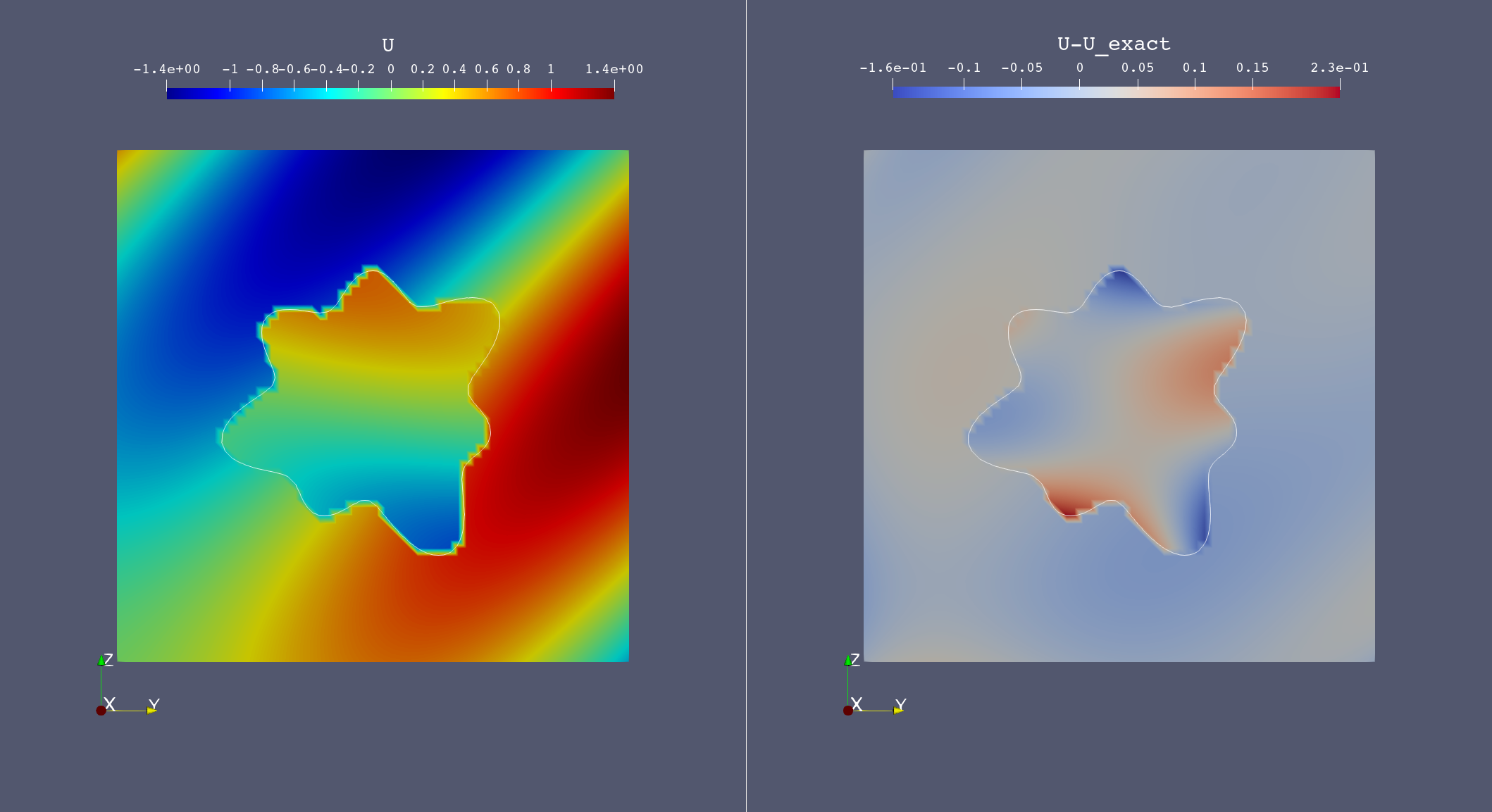}
	\includegraphics[width=\linewidth]{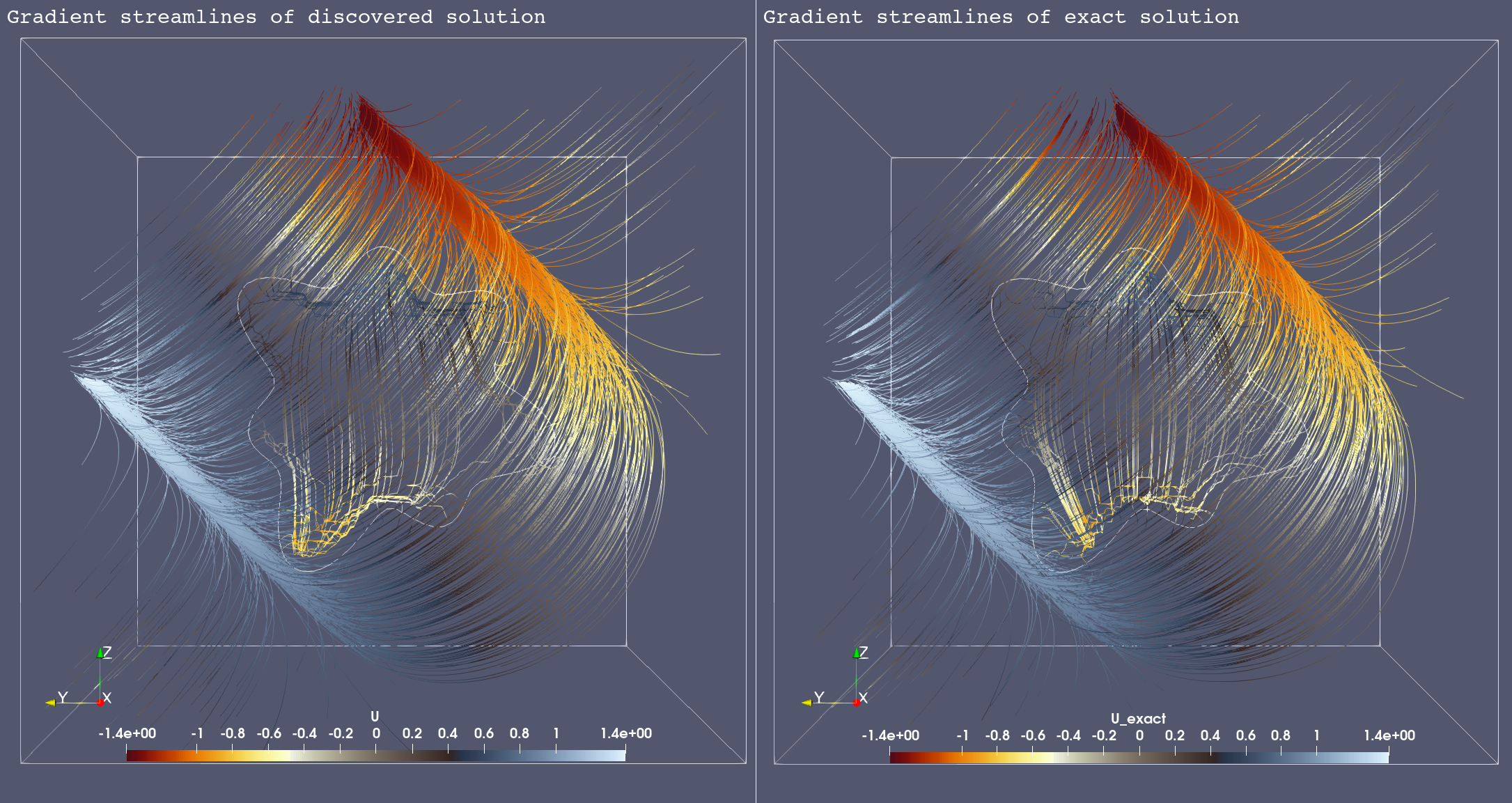}
	\caption{Illustration of exact and numerical solutions (top row) and gradient streamlines (bottom row) on a $64\times 64 \times 64$ grid.}
	\label{fig:star_sol}
\end{figure}

Table \ref{tab:convergence_star} compares approaches I and II for treating the jump conditions using domain switching optimization strategy, these results are illustrated in figure \ref{fig:star}. Table \ref{tab:convergence_star_batching} demonstrates the effect of batching and a multi-resolution training policy on the convergence. We observe convergence is still retained although at a lower rate, and batching can degrade the accuracy at finest resolution (demonstrated by the last row of table \ref{tab:convergence_star_batching} highlighted in red). Although the results generally demonstrate convergence, they are clearly less accurate than what is currently possible by using traditional numerical solvers. In section \ref{sec:future} we describe strategies that may improve these results, importanly using better preconditioners on par with algebraic multigrid methods used in traditional solvers. However, we emphasize that in the case of neural network models considered here the number of degrees of freedom ($\sim 1000$ trainable parameters) is significantly less than the number of degrees of freedom in traditional solvers (for example a grid of $128^3$ has over $2$ million unknown values), besides it is challenging to separate the effect of expressivity of the network (architecure) from the loss composition and the optimization method.

\subsection{Time complexity and parallel scaling on GPU clusters}
We adopt the problem setup presented in \ref{subsec:star}, however with a considerably more challenging geometry of the Dragon presented in \cite{curless1996volumetric}. In this case we used the signed-distance function produced by SDFGen, and initiated an interpolant based on its values.

\begin{figure}
     \centering
     \begin{subfigure}[b]{\textwidth}
         \centering
         \includegraphics[width=\textwidth]{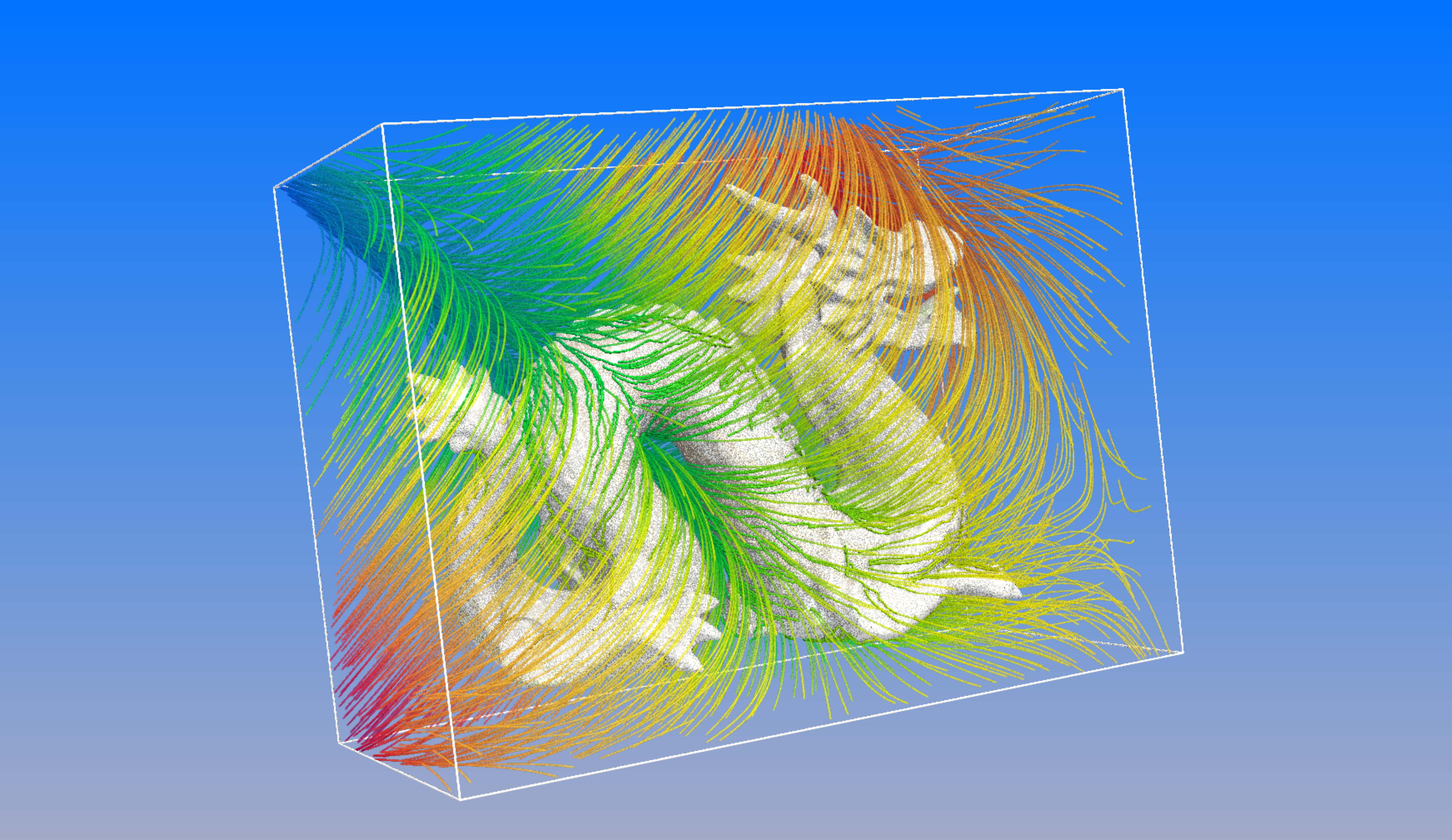}
         \caption{Geometry of the dragon and gradient streamlines, colored by solution values.}
         \label{fig:dragon1}
     \end{subfigure}
     \begin{subfigure}[b]{\textwidth}
         \centering
         \includegraphics[width=\textwidth]{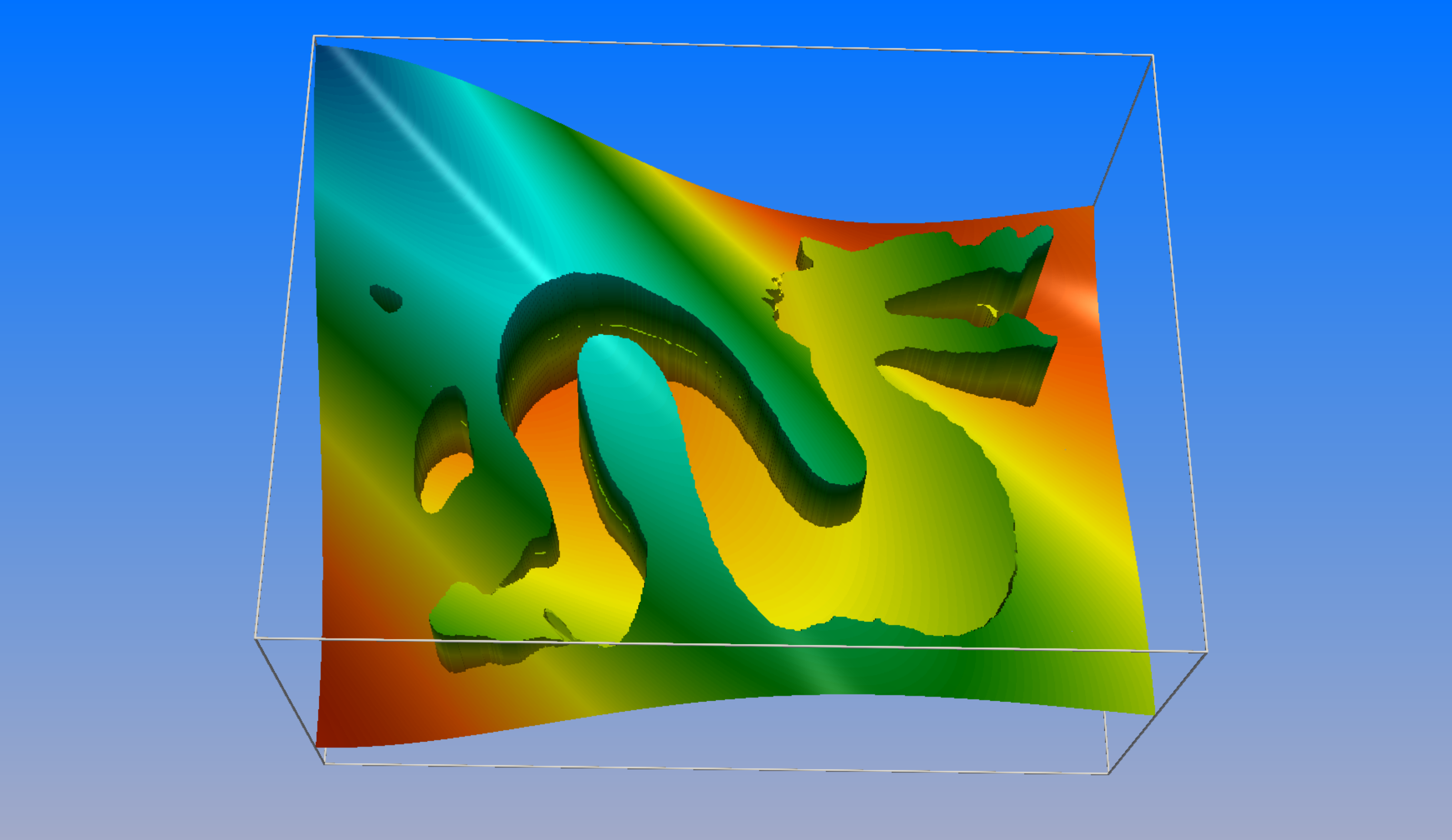}
         \caption{Jump in solution and its gradient are accurately captured by the surrogate neural network model.}
         \label{fig:dragon2}
     \end{subfigure}
        \caption{The NBM approach enables a $1024^3$ effective resolution on a single NVIDIA A6000 GPU. Once trained, it takes sub-milliseconds for the network to evaluate such a simulation that enables near-real-time digital twins for physical systems \cite{tiny-cuda-nn}.}
        \label{fig:dragon}
\end{figure}

The results are shown in figure \ref{fig:dragon}, with a $\rm L^\infty$-error of $0.5$ and RMSE of $0.06$ after 1000 epochs on a base resolution of $64^3$ and implicitly refined onto multi-resolutions $128^3,256^3,512^3$. The neural network pair have only 1 hidden layer with 100 sine-activated neurons, although investigating more complex networks (transformers, symmetry preserving, \textit{etc}.) would likely improve accuracy.

In table \ref{tab:scaling} we report scaling results on NVIDIA A100 GPUs at four base resolutions with three levels of implicit refinement. We used a batchsize of $32\times 32\times 16$ in all cases. At fixed number of GPUs, training time scales linearly (\textit{i.e.}, optimal scaling) with the number of grid points. At a fixed resolution, increasing the number of GPUs accelerates training roughly with $\rm epoch\ time \sim 1/\sqrt{\text{\# }GPUs}$, although the advantage is more effective at higher resolutions. Compile time increases with resolution and decreases with number of GPUs. A maximum grid size of $1024^3$ at multi-resolutions $1024^3,\ 2048^3, 4096^3,\ 8192^3$ was simulated on one NVIDIA DGX with 8 A100 GPUs. The results are shown in figure \ref{fig:dragon}.

\begin{table}[ht]
\begin{center}
\caption{Scaling test. Time per epoch (sec) and JAX compile time for different configurations.} \label{tab:scaling}
\begin{adjustbox}{width=\textwidth}
\begin{tabular}{ccccccccc}
\toprule
base resolution: & \multicolumn{2}{c}{$64^3$}& \multicolumn{2}{c}{$128^3$} & \multicolumn{2}{c}{$256^3$}  & \multicolumn{2}{c}{$512^3$}\\
\cmidrule(r){1-1}\cmidrule(r){2-3}\cmidrule(r){4-5}\cmidrule(r){6-7}\cmidrule(r){8-9}
$\rm A100\ GPUs$ &  epoch  &  compile   &  epoch   & compile &    epoch & compile  &  epoch    &   compile \\
\cmidrule(r){1-1}\cmidrule(r){2-3}\cmidrule(r){4-5}\cmidrule(r){6-7}\cmidrule(r){8-9} 
$1$        & $0.908$ &  $9.027$   & $6.960$  & $9.288$ & $55.287$ & $12.164$ &  $438.45$  &  $49.020$ \\
$2$        & $0.657$ &  $7.575$   & $5.893$  & $7.823$ & $47.360$ & $10.045$ &  $378.98$  &  $39.815$ \\ 
$4$        & $0.405$ &  $7.480$   & $3.629$  & $7.863$ & $28.261$ & $9.129$  &  $226.73$  &  $27.405$ \\ 
$8$        & $0.384$ &  $7.983$   & $3.340$  & $7.901$ & $26.799$ & $9.154$  &  $204.88$  &  $20.632$\\ \bottomrule
\end{tabular}
\end{adjustbox}
\end{center}
\end{table}

\subsection{Comparison with other methods}\label{sec::compare}
Recently other methods have been proposed in the framework of PINNs for solving interfacial PDEs with jump conditions in two \cite{hu2022discontinuity} and three spatial dimensions \cite{wu2022inn}. In this section we compare our method with the recently proposed interfaced neural networks (INN) algorithm of Wu \textit{et al.} (2022) \cite{wu2022inn}. INN uses AD for computing gradients similar to other PINN-like algorithms. INN uses two neural networks for the two computational subdomains and applies an extended multiple-gradient descent (MGD) method in the training phase. This method utilizes information from multiple gradients to adjust and optimize the balance between different terms in the loss function. Intuitively, balancing the weights of different loss terms from different subdomains in the INN method achieves a similar objective as the application of preconditioners in NBM, although NBM preconditions residual contributions from each grid point. Below we compare with example problem 4 of \cite{wu2022inn}.

\subsubsection{Linear Poisson-Boltzmann equation (LPBE) over a sphere}\label{sec:lpbe}
LPBE describes electrostatic potential around biomolecules solvated in low concentration ionic solvents. Here the computational domain is $\Omega\in [-2.5,2.5]^3$ in dimensionless units. The solution contains a regular component and a singular component due to existence of a unit positive charge, $+e_c$, centered at $\mathbf{x}_c=(0,0,0)$. To solve LPBE decomposition techniques are usually employed, here we are only interested in the regular component that constitutes an interfacial problem given by
\begin{align*}
&-\nabla \cdot (\mu^\pm \nabla u) + k^\pm u = 0, & \mathbf{x}\in \Omega^\pm,\\
& [u]=g(\mathbf{x}), &\mathbf{x}\in \Gamma,\\
&[\mu \partial_n u]=\mu^-\partial_n g(\mathbf{x}), &\mathbf{x}\in \Gamma,\\
&u=\frac{\omega}{4\pi\mu^+}\frac{\exp\big(\kappa (\sigma - r)\big)}{(1+\kappa \sigma)r}, & \mathbf{x}\in \partial \Omega, \\
&g(\mathbf{x})=\frac{\omega}{4\pi \mu^-\vert\vert\mathbf{x}-\mathbf{x}_c\vert\vert_2}
\end{align*}
where $\sigma=1$ is the radius of the sphere, $\mu^-=2$ and $\mu^+=80$ while $k^-=0$ and $k^+=\mu^+\kappa^2$ with $\kappa=1.0299\times 10^{-3}$ and $\omega=7.0465\times 10^3$. The exact solution is then given by
\begin{align*}
&u^-=\frac{\omega}{4\pi \sigma}\bigg(\frac{1}{\mu^+(1+\kappa\sigma)} - \frac{1}{\mu^-} \bigg), &u^+=\frac{\omega}{4\pi \mu^+} \frac{\exp\big( \kappa (\sigma - r)\big)}{(1+\kappa\sigma)r}
\end{align*}

Noting that all equations are scaled by $\omega$, we normalize the solution values to the range $[-1.0,1.0]$ by scaling $\omega$ by $\hat{\omega}=293.6$, this is equivalent to setting $\omega=24.0$ in the solver and solving for $\hat{u}^\pm=u^\pm/\hat{\omega}$ instead. This is to avoid the issues with scaling activation functions and having plausible numerical values close to unit value during training. Note that this trick does not impact the \textit{relative} $L_2$ errors for comparison with results reported for INN. For fair comparison with INN we use exactly the same grid point generation code as in INN to generate $\rm N_+=2000$, $\rm N_-=100$, $\rm N_b=1000$, and $\rm N_\Gamma=200$. We build small voxels centered on these grid points with side length of $\rm 0.00244$ corresponding to $\rm 2048$ resolution. We report error estimates in table \ref{tab:convergence_lpbe} using mostly shallow  multilayer perceptron (MLP) architectures as well as the ResNet architecture. 

\subsubsection{Analysis}
Based on table 4 of \cite{wu2022inn}, the smallest network architecture for INN is a pair of one-block ResNet networks with 20 hidden neurons per layer in $\Omega^-$ and 40 neurons in $\Omega^+$, this model produces a relative $\rm L_2$ error of $4.1\times 10^{-2}$ which is on par with the shallow MLP in our experiments, yet NBM trains to this level of accuracy at a much smaller computational cost (see table \ref{tab:convergence_lpbe} for size of our network parameters, timings, and memory requirement). 

On the other hand, tables 4 and 8 of \cite{wu2022inn} report the best relative  $\rm L_2$ value achieved by INN being $3.675\times 10^{-4}$ using a pair of ResNets with 3 blocks each (\textit{i.e.}, each ResNet block contains 2 hidden layers) with 40 neurons per hidden layer to learn the solution in $\Omega^-$ as well as 80 neurons per hidden layer for the solution in $\Omega^+$, with a total of $41,202$ trainable parameters. Regarding validity of these values, we note that the errors reported by \cite{wu2022inn} are estimated on a relatively coarse grid of size $41\times 41\times 41$ that sparsely samples the error in the domain, in contrast our results are measured over a fine resolution of $256\times 256\times 256$ grid points to reliably quantify model accuracy in the space between the points that were used for training. For comparison, we implemented the exact settings described here for the ResNet architecture and used the \texttt{L-BFGS-B} algorithm from the excellent \texttt{jaxopt} package \cite{jaxopt_implicit_diff}, we also removed the $\hat{\omega}$ scaling, \texttt{JAX-DIPS} only takes $89$ seconds to achieve a relative error of $2.22\times 10^{-2}$ using $\rm 3.69\ GB$ of memory. Although the error is still lower than the ones achieved by INN, we emphasize this can improve by using state-of-the-art preconditioners, and also timings prove the efficiency of NBM that is on-par with and even superior to automatic differentiation based algorithms.

\begin{table}[ht]
\begin{center}
\caption{Convergence test of example \ref{sec:lpbe}, positive region uses \texttt{CeLU} activations while negative region uses \texttt{tanh} activations. Note that in all cases there exists an output linear layer. Measurements are on a single NVIDIA A6000 GPU. Errors are measured on a $256\times 256\times 256$ grid to sample the space far from training points. These results are obtained using the Jacobi preconditioner, Adam and L-BFGS-B optimizers, and the discretization of \ref{sec::FD}.} \label{tab:convergence_lpbe}
\begin{adjustbox}{width=\textwidth}
\begin{tabular}{cccccccccc}
\toprule
\multicolumn{3}{c}{model $-\vert+$} & \multicolumn{3}{c}{Accuracy} & \multicolumn{4}{c}{Training Statistics} \\
\cmidrule(r){1-3}\cmidrule(r){4-6}\cmidrule(r){7-10}
type & \#hidden layer:\#hidden unit  & \#params  &   RMSE   &   $\rm L^\infty$ &  Rel. $\rm L^2$  & epochs & time (sec) & Mem. (GB) & optimizer\\
\midrule
\rowcolor{lightgray}\multicolumn{10}{c}{using INN grid with $\rm N_+=2000$, $\rm N_-=100$, $\rm N_b=1000$, {\color{red}$\rm N_\Gamma=200$} }  \\
MLP & $1:1~ \vert 2:10$ & $167$ &  $3.27\times 10^{-3}$ & $1.22\times 10^{-2}$ & $1.92\times 10^{-2}$ & $50,000$ & $693$ & $1.7$ & \texttt{Adam}@\texttt{optax} \\
MLP & $1:1~ \vert 2:10$ & $167$ &  $3.50\times 10^{-3}$ & $1.28\times 10^{-2}$ & $2.15\times 10^{-2}$ & $50,000$ & $95$ & $1.1$ & \texttt{L-BFGS-B}@\texttt{jaxopt} \\ \midrule 
MLP & $1:10 \vert 2:10$ & $212$ &  $3.98\times 10^{-3}$ & $2.13\times 10^{-2}$ & $2.34\times 10^{-2}$ & $50,000$ & $738$ & $1.8$& \texttt{Adam}@\texttt{optax} \\
MLP & $1:10 \vert 2:10$ & $212$ &  $3.43\times 10^{-3}$ & $1.29\times 10^{-2}$ & $2.10\times 10^{-2}$ & $50,000$ & $120$ & $1.17$& \texttt{L-BFGS-B}@\texttt{jaxopt} \\ \midrule
MLP & $2:10~ \vert 2:10$ & $322$ &  $3.83\times 10^{-3}$ & $1.55\times 10^{-2}$ & $2.25\times 10^{-2}$ & $50,000$ & $784$ & $1.9$ & \texttt{Adam}@\texttt{optax}\\
MLP & $2:10~ \vert 2:10$ & $322$ &  $3.52\times 10^{-3}$ & $1.31\times 10^{-2}$ & $2.16\times 10^{-2}$ & $50,000$ & $96$ & $1.2$ & \texttt{L-BFGS-B}@\texttt{jaxopt}\\ 
\rowcolor{lightgray}\multicolumn{10}{c}{using INN grid with $\rm N_+=2000$, $\rm N_-=100$, $\rm N_b=1000$, {\color{red}$\rm N_\Gamma=1000$} }  \\
MLP & $0:1~ \vert 2:10$ & $165$ &  $4.65\times 10^{-3}$ & $1.53\times 10^{-2}$ & $2.73\times 10^{-2}$ & $50,000$ & $923$ & $1.78$ & \texttt{Adam}@\texttt{optax} \\ 
MLP & $0:1~ \vert 2:10$ & $165$ &  $3.52\times 10^{-3}$ & $1.29\times 10^{-2}$ & $2.16\times 10^{-2}$ & $50,000$ & $92$ & $1.15$ & \texttt{L-BFGS-B}@\texttt{jaxopt} \\ \midrule
MLP & $1:1~ \vert 2:10$ & $167$ &  $1.63\times 10^{-3}$ & $4.20\times 10^{-3}$ & $9.62\times 10^{-3}$ & $50,000$ & $946$ & $1.80$ & \texttt{Adam}@\texttt{optax} \\ 
MLP & $1:1~ \vert 2:10$ & $167$ &  $3.53\times 10^{-3}$ & $1.28\times 10^{-3}$ & $2.17\times 10^{-2}$ & $50,000$ & $99$ & $1.18$ & \texttt{L-BFGS-B}@\texttt{jaxopt} \\ \midrule
MLP & $1:10~ \vert 2:10$ & $212$ &  $3.80\times 10^{-3}$ & $1.47\times 10^{-2}$ & $2.22\times 10^{-2}$ & $50,000$ & $977$ & $1.85$  & \texttt{Adam}@\texttt{optax}\\ 
MLP & $1:10~ \vert 2:10$ & $212$ &  $3.51\times 10^{-3}$ & $1.31\times 10^{-2}$ & $2.14\times 10^{-2}$ & $50,000$ & $111$ & $1.18$  & \texttt{L-BFGS-B}@\texttt{jaxopt}\\ \midrule
MLP & $2:10~ \vert 2:10$ & $322$ &  $3.72\times 10^{-3}$ & $1.37\times 10^{-2}$ & $2.19\times 10^{-2}$ & $50,000$ & $1031$ & $1.93$ & \texttt{Adam}@\texttt{optax}\\ 
MLP & $2:10~ \vert 2:10$ & $322$ &  $3.49\times 10^{-3}$ & $1.29\times 10^{-2}$ & $2.14\times 10^{-2}$ & $50,000$ & $126$ & $1.22$ & \texttt{L-BFGS-B}@\texttt{jaxopt}\\  \midrule
MLP & $5:40~ \vert 5:40$ & $13,522$ &  $3.54\times 10^{-3}$ & $1.30\times 10^{-2}$ & $2.17\times 10^{-2}$ & $50,000$ & $80$ & $3.34$ & \texttt{L-BFGS-B}@\texttt{jaxopt}\\  \midrule\midrule
ResNet & $3:80~ \vert 3:40$ & $41,202$ &  $-$ & $-$ & $2.22\times 10^{-2}$ & $50,000$ & $89^\ast$ & $3.69$ &\texttt{L-BFGS-B}@\texttt{jaxopt} \\ 
\rowcolor{lightgray}\multicolumn{10}{c}{uniform $32\times 32\times 32$ grid}  \\
MLP & $1:1 \vert 2:10$ & $167$ & $3.15\times 10^{-3}$  & $1.3\times 10^{-2}$ & $1.85\times 10^{-2}$  & $10,000$ & $794$ & $3.15$ &\texttt{Adam}@\texttt{optax} \\
MLP & $2:10 \vert 2:10$ & $322$ & $2.99\times 10^{-3}$  & $1.25\times 10^{-2}$ & $1.83\times 10^{-2}$  & $50,000$ & $160$ & $2.96$ &\texttt{L-BFGS-B}@\texttt{jaxopt} \\
\bottomrule
\end{tabular}
\end{adjustbox}
\end{center}
\end{table}

\section{Discussion and future work}

\subsection{Spatial gradient calculation using finite discretization improves neural network regularity} 
NBM minimizes the sum of discretization residuals calculated over a set of grid voxels centered on unstructured grid points using any arbitrary finite discretization scheme. Essentially, the finite discretization scheme imposes specific spatial regularities to the neural network predictions inside the spatial volume enclosed by each grid voxel centered at each grid point. This adds an extra \textit{spatial regularization} on the neural network predictions. In contrast to PINNs, computing loss using finite discretizations imposes a spatial regularity on the neural network predictions inside each finite volume; \textit{i.e.}, neural network predictions at vertices and faces of the grid voxel are related to the voxel center according to finite discretization equations. This adds more spatial regularity to the trained neural network compared to automatic differentiation that computes gradients in completely local manner. 

One important advantage of this extended \textit{spatial support} for computing gradients is realized when using conservative finite discretizations such as the one in section \ref{sec::FD}. These explicit conservation laws help to constrain the neural network within each finite volume. We believe this spatial regularization is an important mechanism to encode physical and mathematical priors into the training of neural networks. During training on a set of grid points the neural network not only learns the solution at each point, but also is informed about local symmetries within its neighborhood. As a result of spreading the residual probe according to finite discretization equations, the regularity of the network in the regions between training points follows mathematical regularities and symmetries that are encoded in the bootstrapped finite discretization scheme. In the case of interfacial PDEs with jump conditions considered here, these conservation laws are explicitly considered in section \ref{sec::FD} that govern the solution flux across the interface given the jump conditions. We believe this \textit{mathematically-informed spatial regularization} is responsible for explaining the improved regularity in errors and convergence rates observed in our numerical experiments of section \ref{sec:examples}.

Interestingly, the numerical results in section \ref{sec::compare} demonstrate that even shallow multilayer perceptrons with a few hundreds of trainable parameters can in practice reach levels of accuracy that otherwise require complex architectures such as ResNets with tens of thousands of parameters using PINNs. A network's ability to converge more accurately depends on the optimization strategy, because a simple MLP is a universal function approximator and therefore can in theory represent any function to any level of accuracy. Here we have proposed a new way for optimizing neural networks for PDEs that encodes solutions more efficiently. This is a result of more constrained training that underlies NBM residuals.

\subsection{Training with only first order automatic differentiation improves taining performance}
Additionally, NBM distinguishes between spatial gradients of the neural network with respect to its inputs (needed for the PDE evaluation) from gradients of the neural network with respect to its internal parameters (needed for optimization). NBM computes spatial gradients using finite discretization schemes, while optimization of its internal parameters are handled by the AD-based first order optimizers. This removes the need for higher order AD computations over the computational graph of the neural network during training, therefore increasing scalability and decreasing memory requirements for training deeper and more complex neural network models with higher number of trainable parameters.

The numerical results presented in section \ref{sec:examples} and comparisons of section \ref{sec::compare} demonstrate the superior computational efficiency of NBM implemented in \texttt{JAX-DIPS} compared to PINN-like frameworks such as INN. Our results demonstrate that by removing the expensive higher order AD calculations we gain $\rm \sim 10x$ speedup while maintaining minimal memory requirements ($\rm 1-2\ GB$) for training neural network models.

\subsection{Current shortcomings and future improvements}\label{sec:future}
There are several algorithmic improvements and applications for the current work that we did not explore in this manuscript. The most important algorithmic improvement is application of more advanced preconditioners that can effectively improve the condition number of linear systems arising in interfacial PDE problems such as the ones considered here. Preconditioning the discretization residuals before applying optimization step is an essential requirement for enhancing accuracy of interfacial PDEs in finite discretization methods. The irregular geometries of interfaces as well as the jump conditions and discontinuous parameters of the environment lead to ill-conditioned linear systems. In the current work we only considered the Jacobi preconditioner which is the most basic preconditioner due to ease of implementation. We believe the most important improvement to the current work is application of more advanced preconditioners such as the algebraic multigrid preconditioner \cite{stuben2001introduction,falgout2006introduction}. The numerical methods considered in section \ref{sec::FD} can in principle reach several orders of magnitudes better numerical accuracies when preconditioned properly, see \cite{bochkov2020solving,guittet2015solving} for results of using algebraic preconditioners. These reported results also indicate the obtained convergence results in section \ref{sec:examples} are not simply reflective of the accuracy limitation of the bootstrapped finite discretization method, but they are limited by the improper condition number of the residuals and suboptimal neural network models.

Adaptive mesh refinement is another important algorithmic upgrade to the current work. Adaptive grids with enhanced resolutions closer to the interface and coarser grid cells in the bulk play a significant role in both reducing the computational load as well as improving accuracy of fluxes across interfaces by increasing contributions to the total loss from points closer to the interface. Additionally, more expressive neural network architectures should be considered in \texttt{JAX-DIPS} by adding to the model class of the library. We only considered MLPs, however in recent years there have been a plethora of deep neural network models that have shown great promise such as FNO \cite{li2020fourier}, DeepONet \cite{lu2019deeponet}, and their numerous variants.

In terms of applications, our ultimate goal is training neural operators that can map from different geometries for discontinuities to the solution field. This is critical for developing near real-time simulations of time-evolving systems in digital twins for physical systems. NBM is applicable for physics-driven training of neural operators for elliptic problems with freely moving boundaries, we will present this work in a future time. Additionally, benefiting from differentiablity of the solver we will explore utility of NBM for solving inverse-PDE problems as well as parameterized PDEs.

\section{Conclusion}
\label{sec:conc}
We developed a differentiable multi-GPU framework for solving partial differential equations with jump conditions across irregular interfaces in three spatial dimensions. We developed the neural bootstrapping method to leverage existing finite discretization methods for optimization of neural network internal parameters, while explicitly calculating the spatial gradients using advanced finite discretization methods that encode symmetries and conservation laws governing the PDE solution and its flux across interfaces. Importantly, our framework only uses first order automatic differentiation for optimizing internal state of the neural networks, this technique provides an efficient alternative for training higher order PDE systems by avoiding computational challenges posed by higher order AD over deep neural networks. Moreover, NBM paves the path for obtaining more accurate neural network models of PDEs by leveraging numerical preconditioners by, intuitively, regularizing residuals computed on individual grid points, thus improving optimization gradients.


\section*{Acknowledgement}
This work has been partially funded by ONR N00014-11-1-0027. We extend our sincere gratitude to the reviewers for their invaluable contributions and insightful feedback, which greatly enhanced the quality and rigor of our work.

\newpage
\bibliographystyle{abbrv}
\addcontentsline{toc}{section}{\refname}
\bibliography{references}

\appendix

\section{Geometry Module}
\label{appendix:geom}
The geometric representation of interfaces are handled using the level-set method which is described in subsection \ref{appendinx:lvlset}. The level set methodology provides the foundation for interpolations as well as geometric calculations that are described in subsections \ref{appendix:interp} and \ref{appendix:geometry}, respectively.

\subsection{Level-set method}\label{appendinx:lvlset}
The level-set method for solving free boundary problems was introduced by Osher \& Sethian (1988) \cite{osher1988fronts}. The sharp interface is described as the zero contour of a signed-distance function, $\phi$, whose evolution is given by the advection equation according to some velocity field dictated by the physics of the problem, $\mathbf{v}$, that is defined over the moving boundary
\begin{align*}
	\frac{\partial \phi}{\partial t} + \mathbf{v}\cdot \nabla \phi = 0
\end{align*}
This \textit{implicit} representation of the moving boundaries resolves the need for the challenging task of adapting the underlying grid to the yet-unknown discontinuities in the solution field. The computational simplicity of using Cartesian grids for solving free boundary problems with irregular geometries, as well as the ability to simulate freely moving discontinuities in a \textit{sharp}-manner are the two main advantages of the level-set method for this class of PDE problems. Besides implicit representation of the free boundaries, the level-set function can be used to compute normal vectors to the interface 
\begin{align*}
	\mathbf{n} = \nabla\phi / \vert \nabla \phi \vert 
\end{align*}
as well as the curvature of the interface
\begin{align*}
	\kappa = \nabla \cdot \mathbf{n}
\end{align*}
Details of the interpolation and geometric integration methods in the level-set module are elaborated in \ref{appendix:geom}.

\subsection{Interpolation methods}\label{appendix:interp}
An important building block is the ability to interpolate field values anywhere inside a grid cell given the values on the grid points. In \texttt{JAX-DIPS} we currently support two types of interpolation schemes that have been used in the context of the level set method for achieving second-order accurate solutions by Min \& Gibou (2007a)\cite{MIN2007300}: (i) trilinear interpolation, and (ii) quadratic non-oscillatory interpolation.

\subsubsection{Trilinear interpolation}
In a unit grid cell, rescaled to $\mathcal{C}\in [0,1]^3$, the trilinear interpolation at a point $(x,y,z)\in \mathcal{C}$ uses the grid values on the parent cell vertices according to equation 11 of \cite{MIN2007300},
\begin{align*}
 \phi(x,y,z)= \sum_{i,j,k\in {0,1}} \phi(i,j,k) (-1)^{i+j+k}(1-x-i)(1-y-j)(1-z-k)
\end{align*}
Trilinear interpolation is based on polynomials of order 1 and offers accuracy of order 2 using 8 immediate vertices in a grid cell.

\subsubsection{Quadratic non-oscillatory interpolation}
Quadratic interpolation extends the trilinear interpolation by adding second order derivatives of the interpolating field. This is needed because trilinear interpolation is sensitive to presence of discontinuities and kinks which are abundant in the context of free boundary problems. The extension reads
\begin{align*}
	\phi(x,y,z)= &\sum_{i,j,k\in {0,1}} \phi(i,j,k) (-1)^{i+j+k}(1-x-i)(1-y-j)(1-z-k) \\
				& - \phi_{xx} \frac{x(1-x)}{2} - \phi_{yy} \frac{y(1-y)}{2}- \phi_{zz} \frac{z(1-z)}{2} 
\end{align*}
where second order derviatives are sampled as the minimum value on the parent cell vertices to enhance numerical stability of the interpolation
\begin{align*}
	\phi_{xx} & = \min_{v \in \textrm{vertices}(\mathcal{C})} \vert D_{xx} \phi_v \vert \\
	\phi_{yy} & = \min_{v \in \textrm{vertices}(\mathcal{C})} \vert D_{yy} \phi_v \vert \\
	\phi_{zz} & = \min_{v \in \textrm{vertices}(\mathcal{C})} \vert D_{zz} \phi_v \vert
\end{align*}
The second order derivative operator is the familiar 5-point finite difference stencil.

\subsection{Geometric integration}\label{appendix:geometry}
\subsubsection{Integration over 3D surfaces and volumes}
We use uniform Cartesian grids. For computational cells that are crossed by the interface, \textit{i.e.} $\mathcal{V}_{i,j,k}\cap \Gamma \neq 0$, we use the geometric integrations proposed by Min \& Gibou (2007b) \cite{min2007geometric}. In this scheme each grid cell, $\mathcal{C}$, is decomposed into five tetrahedra by the middle-cut triangulation \cite{sallee1984middle} (each cell is rescaled to $[0,1]^3$) that are described below:
\begin{align*}
	 & \rm S_1 \equiv conv(P_{000} ; P_{100} ; P_{010} ; P_{001}) & \rm{ x = 0\ face,\ y = 0\ face,\ z = 0\ face} \\
	 & \rm S_2 \equiv conv(P_{110} ; P_{100} ; P_{010} ; P_{111}) & \rm{ x = 1\ face,\ y = 1\ face,\ z = 0\ face} \\
	 & \rm S_3 \equiv conv(P_{101} ; P_{100} ; P_{111} ; P_{001}) & \rm{ x = 1\ face,\ y = 0\ face,\ z = 1\ face} \\
	 & \rm S_4 \equiv conv(P_{011} ; P_{111} ; P_{010} ; P_{001}) & \rm{ x = 0\ face,\ y = 1\ face,\ z = 1\ face} \\
	 & \rm S_5\equiv conv(P_{111} ; P_{100} ; P_{010} ; P_{001})  & \rm{ no\ face\ exposure}
\end{align*}

Hence each 3D grid cell is the union of $5$ tetrahedra (simplices) $\mathcal{C}=\cup_{i=1}^5 S_i$, where each simplex is identified by the pre-existing vertices of the grid cell (hence not creating new grid points). Given the values of the level set function sampled at these vertices one can compute coordinates of intersection points of the interface with each of the simplices $S_i \cap \Gamma $ as well as the negative domain $S_i \cap \Omega^-$. If ${P_0,\cdots, P_3}$ are the four vertices of a simplex $S$, then $\Gamma$ crosses an edge $P_i P_j$ if and only if $\phi(P_i)\phi(P_j)<0$ and the intersection point across this edge is given by linear interpolation:
\begin{align*}
P_{ij}=P_j \frac{\phi(P_i)}{\phi(P_i) - \phi(P_j)} - P_i \frac{\phi(P_j)}{\phi(P_i) - \phi(P_j)}
\end{align*}

Number of negative level-set values on the $4$ (in 3D) tetrahedron vertices classifies the specific configuration for intersection between simplex $S$ and the interface through a variable $\eta(\phi, S)=n({P_i | \phi(P_i)<0})$. In 3D, possible values are $\eta =0,1,2,3,4$ that correspond to the four configurations for the intersection cross section enumerated below:
\begin{itemize}
\item $S\cap\Gamma$:
	\begin{itemize}
		\item $\eta=0$: tetrahedron ($S$) is completely in positive domain with no intersection, $S\cap\Gamma=\emptyset$.
		\item $\eta=1$: with a single vertex in negative domain and remaining three in positive domain, the tetrahedron and interface have exactly $3$ intersection points, the simplex $S\cap\Gamma$ has exactly $3$ vertices; \textit{cf.}, see figure 2 (center) of  \cite{min2007geometric}.   
		\item $\eta=2$: with two vertices in negative domain and remaining two in positive domain, the cross section has four vertices that is splitted into two 3-vertex simplices; \textit{cf.}, see figure 2 (right) of  \cite{min2007geometric}.
		\item $\eta=3$: with one vertex in positive domain and remaining three vertices in negative domain, the cross section has 3 vertices that is computed by inverting the sign of the level-set values on vertices and following the instruction for case $\eta=1$.
		\item $\eta=4$: tetrahedron is completely in negative domain with no intersection, $S\cap\Gamma=\emptyset$.
	\end{itemize}

\item $S\cap \Omega^-$:
	\begin{itemize}
		\item $\eta=0$: tetrahedron is completely in positive domain with no intersection, $S\cap\Omega^-=\emptyset$; 
		\item $\eta=1$: the intersection $S\cap\Omega^-$ is characterized by a single tetrahedron with 4 vertices according to figure 4 (left) of \cite{min2007geometric}; \textit{i.e.}, one vertex is the negative level-set vertex of the parent tetrahedron and three others are interpolated points over the three edges pertaining to the negative vertex.
		\item $\eta=2$: the intersection $S\cap\Omega^-$ is characterized by three tetrahedra with 12 vertices according to figure 4 (center) of \cite{min2007geometric}. Note that there 
		\item $\eta=3$: the intersection $S\cap\Omega^-$ is characterized by three tetrahedra with 12 vertices according to figure 4 (right) of \cite{min2007geometric}.
		\item $\eta=4$: tetrahedron is completely in negative domain and $S\cap\Omega^-=S$; 
	\end{itemize}

\end{itemize}

Note that although we only need to allocate memory for at most 4 vertices to uniquely identify $S\cap\Gamma$, in \texttt{JAX-DIPS} we choose to pre-allocate memory for two 3-vertex simplex data structure per $S$ with a total of $6$ vertices to separately store information for the cross section geometry. Similarly for $S\cap \Omega^-$ we pre-allocate memory for a three 4-vertex simplex data structure per $S$. Altogether, in the current implementation the geometric information of intersection points for each simplex $S$ is expressed in terms of 5 simplicies (2 three-vertex simplices for surface area and 3 four-vertex simplicies for volume) using 18 points; this is an area for future optimization.

Having the intersection points, we compute surface and volume integrals of a given field over the interface $\Gamma$ and in negative domain $\Omega^-$ as a summation of integrals over the identified simplices. For each simplex (with $n=3$ or $n=4$ vertices) surface and volume integrals can be numerically computed by having these vertices ${P_0,\cdots, P_n}$ and the values of the field $f$ at these vertices according to 
\begin{align*}
\int_S f dx = \textrm{vol}(S)\cdot \frac{f(P_0) + \cdots + f(P_n)}{n+1}
\end{align*}
where
\begin{align*}
\textrm{vol}(S) = \frac{1}{n!} \bigg\vert \textrm{det} \begin{pmatrix}
(P_1-P_0)\hat{e}_1 & \cdots & (P_n-P_0)\hat{e}_1 \\
\vdots  & \ddots & \vdots  \\
(P_1-P_0)\hat{e}_n  & \cdots & (P_n-P_0)\hat{e}_n 
\end{pmatrix} \bigg\vert
\end{align*}
with $\hat{e}_i$ being the $i^{th}$ Cartesian unit basis vector.

\subsubsection{Cross sections of interface with grid cell faces}
For the numerical discretizations considered in this work we also need the surface areas for simplices created at the intersection of $\Gamma$ with each of the 6 faces of a grid cell $\mathcal{C}$. In \texttt{JAX-DIPS} for each face we reuse two corresponding simplices exposed to that face that were calculated in the geometric integrations module, explicitly: 
\begin{itemize}
	\item $x=0$ face has contributions from ($S_1$, $S_4$)
	\item $x=1$ face has contributions from ($S_2$, $S_3$)
	\item $y=0$ face has contributions from ($S_1$, $S_3$)
	\item $y=1$ face has contributions from ($S_2$, $S_4$)
	\item $z=0$ face has contributions from ($S_1$, $S_2$)
	\item $z=1$ face has contributions from ($S_3$, $S_4$)
\end{itemize}
For each face, we extract vertices from $(S_i,S_j)$-pair that lie on the considered face and sum the surface areas in (negative domain) contributed from each simplex on that face; therefore, the portion of the face surface area in the positive domain is simply the complementing value $\rm area^+=area_{face}-area^-$; \textit{i.e.}, this ensures sum of areas adds up to the exact face area in downstream computations.

\section{Solving interface problems with physics-informed neural networks}

Physics-informed neural networks (PINNs) are based on minimizing a penalty function on evaluated on a point cloud, where the actual PDE constitutes the residual and automatic differentiation is used to compute the spatial derivaties in the PDE. Each of the two relations \ref{eq::extrapolate1}--\ref{eq::extrapolate2} offers an extra residual term that can be used to penalize the loss function in PINNs
\begin{align}
\mathcal{L}_{interface} &= \bigg|\bigg| u^+_{i,j} - u^-_{i,j} - \alpha(\mathbf{r}_{i,j}^{pr}) - \delta_{i,j}\bigg(\big(\frac{\mu^-}{\mu^+}-1\big)\partial_\mathbf{n}u^-(\mathbf{r}^{pr}_{i,j}) + \frac{\beta(\mathbf{r}_{i,j}^{pr})}{\mu^+} \bigg) \bigg|\bigg|_2^2 \\
\intertext{or}
\mathcal{L}_{interface} &= \bigg|\bigg| u^+_{i,j} - u^-_{i,j} - \alpha(\mathbf{r}_{i,j}^{pr}) - \delta_{i,j}\bigg(\big(1-\frac{\mu^+}{\mu^-}\big)\partial_\mathbf{n}u^+(\mathbf{r}^{pr}_{i,j}) + \frac{\beta(\mathbf{r}_{i,j}^{pr})}{\mu^-} \bigg) \bigg|\bigg|_2^2
 \label{eq::interfaceLoss}
\end{align}

Far from interface, the usual procedure for physics-informed neural networks is applicable, namely,
\begin{align*}
\mathcal{L}_{bulk} &= \bigg|\bigg| k^{\pm}u^{\pm} - \nabla \cdot (\mu^{\pm}\nabla u^\pm) - f^{\pm} \bigg|\bigg|_2^2 \\
\mathcal{L}_{boundary} &= \bigg|\bigg|u(\mathbf{r}_{bc}) - \hat{u}(\mathbf{r}_{bc})\bigg|\bigg|_2^2
\end{align*}
hence, the overall loss function for this class of problems shall be
\begin{align*}
\mathcal{L} = \mathcal{L}_{bulk} + \mathcal{L}_{boundary} + \mathcal{L}_{interface}
\end{align*}

\end{document}